\documentclass[11pt,reqno]{amsart}

\usepackage{amsmath}
\usepackage{amsthm}
\usepackage{hyperref}
\usepackage[margin=1.25in]{geometry}

\usepackage{color}

\numberwithin{equation}{section}

\newtheorem{prop}{Proposition}
\newtheorem{lemma}[prop]{Lemma}
\newtheorem{thm}[prop]{Theorem}
\newtheorem{cor}[prop]{Corollary}

\numberwithin{prop}{section}

\theoremstyle{definition}
\newtheorem{defn}[prop]{Definition}

\newtheorem{rmk}[prop]{Remark}

\newcommand{\dt}{\frac{\partial}{\partial t}}
\newcommand{\brs}[1]{\left| #1 \right|}
\newcommand{\gG}{\Gamma}
\newcommand{\gD}{\Delta}
\newcommand{\gd}{\delta}
\newcommand{\gs}{\sigma}
\newcommand{\gw}{\omega}
\newcommand{\ga}{\alpha}
\newcommand{\gb}{\beta}
\newcommand{\gl}{\lambda}
\renewcommand{\ge}{\epsilon}
\newcommand{\N}{\nabla}
\renewcommand{\bar}[1]{\overline{#1}}
\newcommand{\del}{\partial}

\newcommand{\bj}{\bar{j}}
\newcommand{\bk}{\bar{k}}
\newcommand{\bl}{\bar{l}}
\newcommand{\bm}{\bar{m}}
\newcommand{\bn}{\bar{n}}
\newcommand{\bp}{\bar{p}}

\newcommand{\til}[1]{\widetilde{#1}}

\newcommand{\PP}{\mathcal K}

\newcommand{\ohat}[1]{\overset{\circ}{#1}}

\DeclareMathOperator{\Sym}{Sym}
\DeclareMathOperator{\Rc}{Rc}
\DeclareMathOperator{\Ric}{Ric}
\DeclareMathOperator{\Rm}{Rm}

\DeclareMathOperator{\tr}{tr}

\DeclareMathOperator{\Id}{Id}
\DeclareMathOperator{\id}{id}

\DeclareMathOperator{\End}{End}

\begin{document}

\title[Symplectic Curvature Flow]{Symplectic Curvature Flow}

\begin{abstract} We introduce a parabolic flow of almost K\"ahler structures,
providing an approach to constructing canonical geometric structures on
symplectic manifolds.  
We exhibit this flow as one of a
family of parabolic flows of almost Hermitian structures, generalizing our
previous work on parabolic flows of Hermitian metrics.  We exhibit a long time
existence obstruction for solutions to this flow by showing certain smoothing
estimates for the curvature and torsion.  We end with a discussion
of the limiting objects as well as some open problems related to the symplectic
curvature flow.
\end{abstract}

\author{Jeffrey Streets}
\address{Rowland Hall\\
         University of California, Irvine\\
         Irvine, CA 92617}
\email{\href{mailto:jstreets@uci.edu}{jstreets@uci.edu}}

\author{Gang Tian}
\address{BICMR and SMS, Beijing University\\
China and Fine Hall\\
	 Princeton University\\
	 Princeton, NJ 08544}
\email{\href{mailto:tian@math.princeton.edu}{tian@math.princeton.edu}}

\date{August 30th, 2011}

\maketitle

%%%%%%%%%%%%%%%%%%%%%%
\section{Introduction}
%%%%%%%%%%%%%%%%%%%%%%

In the past two decades the study of symplectic manifolds has been very active.
New tools have been introduced
and new insight has been provided, for instance, 
by Gromov's work on pseudoholomorphic curves and symplectic topology
\cite{Gromov}, the Gromov-Witten invariants and their applications to
mirror symmetry (see for instance \cite{ruantian} etc.), Taubes' works on the
Seiberg-Witten equations on symplectic manifolds \cite{Taubes1}, \cite{Taubes2},
invariants coming
from studying Hamiltonian dynamics and Lagrangian intersections (see for
instance
\cite{Hoferbook} and \cite{fooo}).  These approaches have all had a profound
impact on our
understanding of symplectic manifolds, and are linked in the sense that they are
all ``topological'' in nature. The purpose of this paper is to introduce a
geometric approach to studying symplectic manifolds.
Specifically we introduce a new curvature flow which preserves symplectic
structures and evolves almost K\"ahler structures, which always exist on
symplectic manifolds, towards certain canonical geometric structures on
symplectic manifolds.  Hopefully, this curvature flow provides us a very
different approach to and enables
us to apply the methods of geometric analysis to understanding the topology
and geometry of symplectic manifolds from a different point of view.

To begin, let $(M^{2n}, \omega)$ denote a
compact smooth manifold with closed, nondegenerate $2$-form $\omega$.  Any such
$\omega$ admits compatible almost complex structures.  Below we will define a
coupled degnerate parabolic system of equations for a compatible pair $(\omega,
J)$ preserving the symplectic condition for $\omega$.  If the initial almost
complex structure is in fact integrable, then the resulting one-parameter family
of complex structures is fixed, i.e. $J(t) = J(0)$, and the family of K\"ahler
forms $\omega(t)$ is a solution to K\"ahler Ricci flow. This parabolic system is
furthermore a special instance of a general family of parabolic flows of almost
Hermitian structures.  We begin by describing this more general setup, then
proceed to define the flow of almost K\"ahler structures.

Let $(M^{2n}, \gw, J)$ be an almost Hermitian manifold.  Let $\N$ denote the
Chern connection
associated to $(\gw, J)$, which is the unique connection satisfying
\begin{align*}
\N \gw \equiv 0, \qquad \N J \equiv 0, \qquad T^{1,1} \equiv 0
\end{align*}
where $T^{1,1}$ refers to the $(1,1)$ component of the torsion of $\N$ thought
of as a section of $\Lambda^2 \otimes TM$.  Let $\Omega$ denote the
$(4,0)$-curvature tensor associated to this connection, and let
\begin{align*}
S_{ij} = \omega^{kl} \Omega_{klij}.
\end{align*}
Furthermore, let $Q$ denote a $(1,1)$ form which is a quadratic expression in
the torsion $T$ of $\N$.  Let
\begin{align*}
\mathcal K_j^i = \omega^{kl} \N_k N_{l j}^i.
\end{align*}
where $N$ denotes the Nijenhuis tensor associated to $J$.  Also, let $\mathcal
H$ denote a generic quadratic expression in the Nijenhuis tensor which is an
endomorphism
of the tangent bundle which skew-commutes with $J$.  Finally, let
\begin{align} \label{Hdef1}
H = \frac{1}{2} \left[ \gw(\mathcal K - \mathcal H, J) + \gw(J, \mathcal K -
\mathcal H) \right].
\end{align}
These definitions are spelled out in greater detail in the rest of the paper.
Consider the
initial
value problem
\begin{gather} \label{almostcomplexflow}
\begin{split}
\dt \gw =&\ - S + Q + H\\
\dt J =&\ - \mathcal K + \mathcal H\\
\omega(0) =&\ \omega_0\\
J(0) =&\ J_0.
\end{split}
\end{gather}
\noindent This is a
degenerate parabolic system of equations for $(\omega, J)$, with
degeneracy arising from the action of the diffeomorphism group.  In section 3 we
prove the general short-time existence of solutions of
(\ref{almostcomplexflow}), a generalization of Theorem 1.1 of \cite{ST1}.

\begin{thm} \label{almosthermitianste} Let $(M^{2n}, \gw_0, J_0)$ be a compact
almost Hermitian manifold.  There
exists $\ge > 0$ and a unique one parameter family of almost Hermitian
structures $(\omega(t), J(t))$ solving (\ref{almostcomplexflow}) with
initial
condition
$(\gw_0, J_0)$.  If $J_0$ is integrable, then $J(t) = J_0$
for all $t \in
[0, \ge)$.  Furthermore, if $J_0$ is integrable and $g_0$ is K\"ahler, then
$g(t)$ is K\"ahler for all $t \in [0, \ge)$ and $g(t)$ solves the K\"ahler-Ricci
flow with initial condition $g_0$.
\end{thm}

\begin{rmk} It is important to note that equation (\ref{almostcomplexflow}) is
defining a \emph{family} of equations.  Indeed, the choice of $Q$ and $\mathcal
H$ are \emph{arbitrary} in the definition of (\ref{almostcomplexflow}) and the
proof of Theorem \ref{almosthermitianste}.
\end{rmk}

\begin{rmk} When $J_0$ is integrable, the one-parameter family of metrics
$\omega(t)$ is a solution to \emph{Hermitian curvature flow}, as defined in
\cite{ST1}.  Again, the torsion term $Q$ can be arbitrary for the result of
Theorem \ref{almosthermitianste}.
\end{rmk}

\begin{rmk} As will be clear from Proposition \ref{goperator}, it is possible to
define a parabolic flow of metrics compatible with any given almost complex
structure.  Specifically, given $(M^{2n}, J)$ an almost complex manifold, one
can set
\begin{align} \label{almflow}
\dt \omega =&\ - S + Q
\end{align}
where again $Q$ is a $(1,1)$ form which is a quadratic expression in the
torsion.  This viewpoint was considered recently by Vezzoni \cite{Ve}.  When $J$
is integrable, this is precisely the family of equations introduced in
\cite{ST1}.  If one is interested in understanding metrics compatible with a
given almost complex structure, (\ref{almflow}) could be a useful tool.
\end{rmk}

We now proceed to define the flow of almost K\"ahler structures.

\begin{defn} An almost Hermitian manifold
$(M^{2n}, \omega, J)$ is \emph{almost K\"ahler} if
\begin{align*}
d \omega = 0.
\end{align*}
\end{defn}

\noindent This
condition is a very natural extension of K\"ahler geometry, and one may
consult \cite{Apostolov} for a nice fairly recent survey of results on these
structures.  Due to their connection with symplectic geometry, almost K\"ahler
structures have become a central area of mathematics (see for instance
\cite{Donaldson}, \cite{LeBrun}).

An almost K\"ahler structure has an associated Levi Civita connection $D$, as
well as a canonical Hermitian connection $\N$
(which
coincides with the Chern connection) with curvature $\Omega$.  Furthermore, one
can define
\begin{align*}
P_{ij} = \omega^{kl} \Omega_{ijkl}.
\end{align*}
By Chern-Weil theory we know that $P \in \pi c_1(M, J)$, and moreover $d P = 0$.
 In analogy with K\"ahler Ricci flow, it is natural to expect that $P$ is the
right operator by which to flow a symplectic structure.  However, $P \notin
\Lambda^{1,1}$, therefore one is forced to attach a flow of $J$ as well to
preserve compatibility of the pair.  Set
\begin{align*}
\mathcal N_i^j =&\ g^{jk} g_{mn} g^{pq} D_p J_r^m J_i^r D_q J_k^n,\\
\mathcal R_i^j =&\ J_i^k \Rc_k^j - \Rc_i^k J_k^j,
\end{align*}
\noindent and consider the initial value problem
\begin{gather} \label{AKflow}
\begin{split}
\dt \omega =&\ - { P}\\
\dt J =&\ - D^* D J \ { + \mathcal N} + \mathcal R\\
\omega(0) =&\ \omega_0\\
J(0) =&\ J_0.
\end{split}
\end{gather}

\begin{thm} \label{AKflowthm} Let $(M^{2n}, \omega_0, J_0)$ be a compact almost
K\"ahler manifold.  There exists $\ge > 0$ and a unique one-parameter family of
almost K\"ahler structures $(\omega(t), J(t))$ solving (\ref{AKflow}) for $t \in
[0, \ge)$.  Moreover, the pair $(\omega(t), J(t))$ is a solution to an equation
of the type
(\ref{almostcomplexflow}), for appropriate choices of $Q$ and $\mathcal H$, with
$H$ defined by (\ref{Hdef1}).  In particular, this instance of equation
(\ref{almostcomplexflow}) preserves the almost K\"ahler condition.  Finally, if
$J_0$ is integrable, $J(t) = J(0)$ for all $t$ and $\omega(t)$ is a solution to
K\"ahler
Ricci flow.
\end{thm}

\begin{rmk} In \cite{LeWang} a certain geometric evolution equation was studied
on symplectic manifolds.  There the perspective taken is that the symplectic
structure $\omega$ is fixed, and then one studies the gradient flow of the
functional of compatible almost complex structures
\begin{align*}
\mathcal F(J) := \int_M \brs{D J}^2 dV
\end{align*}
where the metric defining the quantities above is that associated to $J$ via
$\omega$.  The proof of short time existence of this flow is already
technical, due to certain local obstructions in prescribing the skew-symmetric
part of the Ricci tensor.  Our approach here is different, as we allow both
$\omega$ and $J$ to change.  This seems to have certain advantages, since for
instance the diffeomorphism action is the only obstruction to parabolicity.
Furthermore, our flow is a natural generalization of K\"ahler Ricci flow,
whereas any K\"ahler metric is already a fixed point for this flow.
\end{rmk}

We are able to derive equations for the evolution of curvature and torsion under
solutions of (\ref{almostcomplexflow}) and (\ref{AKflow}).  The general theory
is similar to the case of Hermitian curvature flow, where one requires bounds on
the curvature, torsion, and first derivative of torsion to conclude long
time existence of the flow.  This result is obtained by proving smoothing
estimates for higher derivatives which hold in the presence of these bounds.
For a technical reason explained in section \ref{smoothing}, one is forced to
get $L^2$ smoothing estimates.  Incidentally, this technical problem does not
occur for (\ref{AKflow}), and one obtains the usual pointwise smoothing
estimates (see Theorem \ref{AKsmoothingthm}).

\begin{thm} \label{generalsmoothingthm} Given $m > 0$, there exists $C = C(m,
n)$ such that if $(M^{2n}, \omega(t), J(t))$ is a solution to
(\ref{almostcomplexflow}) on $\left[0, \frac{\ga}{K} \right]$ satisfying
\begin{align*}
\sup_{M \times \left[0, \frac{\ga}{K} \right]} \{ \brs{\Rm}, \brs{T}^2, \brs{D
T}
\} \leq K,
\end{align*}
then
\begin{align*}
\sup_{M \times \left[0, \frac{\ga}{K} \right]} \{ \brs{\brs{ D^m \Rm}}_{L^2}^2,
\brs{\brs{D^{m+1} T}}_{L^2}^2 \} \leq \frac{C K}{t^{\frac{m}{2}}}.
\end{align*}
\end{thm}

\noindent Using these we obtain the long time existence obstruction.

\begin{thm} \label{lteobs} Let $(M^{2n}, \omega_0, J_0)$ be an almost Hermitian
manifold.  There is a unique solution to (\ref{almostcomplexflow}) on a maximal
time interval $[0, \tau)$.  Furthermore, if $\tau < \infty$ then
\begin{align*}
\limsup_{t \to \tau} \{ \brs{\Rm}_{C^0}, \brs{D T}_{C^0}, \brs{T}^2_{C^0} \} =
\infty.
\end{align*}
\end{thm}

Furthermore, one can improve this regularity requirement in the case of
symplectic curvature flow.  This is because of an
a priori estimate
for $\brs{DJ}^2$ which holds when the curvature is bounded.

\begin{thm} \label{symplte} Let $(M^{2n}, \omega_0, J_0)$ be an almost K\"ahler
manifold.
There is a unique solution to (\ref{AKflow}) on a maximal
time interval $[0, \tau)$.  Furthermore, if $\tau < \infty$ then
\begin{align*}
\limsup_{t \to \tau} \brs{\Rm}_{C^0} = \infty.
\end{align*}
\end{thm}

Here is an outline of the rest of the paper.  In $\S$ \ref{bckgrnd} we review
some
basic aspects of almost Hermitian geometry, and recall the Chern connection.  We
recall and generalize some known curvature identities in section $\S$
\ref{identities}.
We give basic calculations on variations of
almost Hermitian structures in $\S$ \ref{variations}.  In $\S$ \ref{generalste}
and $\S$ \ref{symplecticste} we give the proofs of Theorems
\ref{almosthermitianste} and \ref{AKflowthm}.  Evolution equations for the
curvature tensor, torsion, and their derivatives are shown in $\S$
\ref{curvatureev}, and we use these to prove smoothing estimates which are used
to prove Theorems \ref{generalsmoothingthm},
\ref{lteobs} and \ref{symplte} in $\S$ \ref{smoothing}.  In $\S$
\ref{critstruct}, we give a
discussion of some special properties of the limiting metrics of (\ref{AKflow}).
We end in $\S$ \ref{conc} by posing a number of problems related to symplectic
curvature flow.

\textbf{Acknowledgements:} The first author would like to thank Graham Cox,
Zoltan Szabo, Mohommad Tehrani, and Guangbo Xu for several interesting
conversations on this topic.  The authors would also like to thank Vestislav
Apostolov and Tedi Draghici for pointing out connections to their work, and
Heather Macbeth for finding errors in a previous version of this manuscript.

%%%%%%%%%%%%%%%%%%%%%%%%%%%%%%%%%%%%%%%%%%%%%%%%%
\section{Background on Almost Hermitian Geometry} \label{bckgrnd}
%%%%%%%%%%%%%%%%%%%%%%%%%%%%%%%%%%%%%%%%%%%%%%%%%

In this section we review some basic material about almost Hermitian geometry
and various associated connections.  Let $(M^{2n}, J)$ be an almost complex
manifold.
This means that $J$ is an endomorphism of $TM$ satisfying
\begin{align*}
J^2 = - \Id.
\end{align*}
By the theorem of Newlander-Nirenberg \cite{NN}, the almost complex
structure $J$ is \emph{integrable}, i.e. one can find local complex coordinates
at each point, if and only if the Nijenhuis tensor vanishes.  The Nijenhuis
tensor is
\begin{align} \label{Nijenhuis}
N_J (X, Y) =&\ [J X, J Y] - [X, Y] - J[J X, Y] - J[X, J Y].
\end{align}
As in the case of complex manifold, the almost-complex structure $J$ induces
a decomposition of the space of differential forms on $M$ via the eigenspace
decomposition on $TM$.  In particular we
will write
\begin{gather*}
\Lambda^r(M) \otimes \mathbb C = \bigoplus_{p + q = r} \Lambda^{p,q}.
\end{gather*}
Also, for a general two-tensor $W \in T^* M \otimes T^* M$, let
\begin{align*}
W^J(X, Y) =&\ \frac{1}{2} \left(W(X, Y) + W(JX, JY) \right),\\
W^{-J}(X,Y) =&\ \frac{1}{2} \left(W(X, Y) - W(JX, JY) \right)
\end{align*}
denote the projections of $W$ onto $J$-symmetric and $J$-antisymmetric tensors.

The operator $d$ acts on $\Lambda^r$, but in general one does not have $d
\Lambda^{p,q} \subset \Lambda^{p+1, q} \oplus \Lambda^{p,q+1}$, due to the
potential lack of integrability of $J$.  Finally, we will use the operator
\begin{align*}
 d^c : \Lambda^r &\ \rightarrow \Lambda^{r+1}\\
\psi&\ \rightarrow -J d \psi
\end{align*}
where for a differential $r$-form $\phi$ one has
\begin{align*}
 \left(J \phi \right)(X_1, \dots, X_r) =&\ \phi(J X_1, \dots, J X_r).
\end{align*}

Moving to the metric geometry, let $g$ be an almost Hermitian metric on
$M$,
i.e. $g$ satisfies
\begin{align*}
g( \cdot, \cdot) = g(J \cdot, J \cdot).
\end{align*}
Associated to this pair is the K\"ahler form
\begin{align*}
\omega(\cdot, \cdot) = g(J \cdot, \cdot).
\end{align*}

Next we consider connections associated to almost Hermitian manifolds.  A very
thorough discussion of these connections can be found in \cite{GauduchonConn}.
A linear
connection $\N$ on $TM$ is called \emph{Hermitian} if
\begin{align*}
\N \gw \equiv 0, \qquad \N J \equiv 0.
\end{align*}
These two conditions alone do not suffice to determine a unique connection in
general.  Indeed, there is freedom yet of $\psi \in \Lambda^3(\mathbb R) \cap
\Lambda^{2,1} \oplus \Lambda^{1,2}$ and $B \in \Lambda^{1,1} \otimes TM$
satisfying a certain Bianchi identity (see \cite{GauduchonConn} Proposition 2).
Certain members of this family are chosen according to certain desirable
properties of the torsion.  Frequently, one chooses the Chern connection.

\begin{defn} \label{Cherndef} Given $(M^{2n},\omega, J)$ an almost-Hermitian
manifold, the \emph{Chern
connection associated to $(\omega, J)$} is the unique connection $\N$ satisfying
\begin{align*}
\N \omega \equiv&\ 0\\
\N J \equiv&\ 0\\
T^{1,1} \equiv&\ 0
\end{align*}
where $T$ denotes the torsion tensor of $\N$ and $T^{1,1}$ is the projection of
the vector-valued torsion two-form onto the space of $(1,1)$-forms.
\end{defn}

Gauduchon \cite{GauduchonConn} has identified a canonical family of Hermitian
connections associated to an almost Hermitian pair.   Before describing it
though let
us introduce a further piece of notation.  For $\phi \in \Lambda^{3}$, let
\begin{align*}
\phi^+ :=&\ \phi^{(2,1) + (1,2)}\\
\phi^- :=&\ \phi^{(3,0) + (0,3)}.
\end{align*}
By making certain natural assumptions about the torsion of a Hermitian
connection (see \cite{GauduchonConn} Definition 2), one identifies a
one-parameter family of such connections (\cite{GauduchonConn} 2.5.4)
\begin{gather} \label{generalconnection}
 \begin{split}
\left< \N_X Y, Z \right> =&\ \left< D_X Y, Z \right> + \frac{1}{2} \left< (D_X
J) J Y, Z \right> + \frac{t}{4} \left( \left( d^c \omega \right)^+_{X,Y,Z} +
\left(
d^c \omega \right)^+_{X,JY,JZ} \right).
 \end{split}
\end{gather}
In the formula above $D$ denotes the Levi-Civita connection.  The choice $t = 1$
corresponds to the \emph{Chern connection}.  As a final
important remark we observe that in the case of almost K\"ahler manifolds, this
family reduces to a single point, i.e. there is a
\emph{canonical}
Hermitian connection on almost K\"ahler manifolds, taking the simple form
\begin{gather} \label{AKconnection}
\left< \N_X Y, Z \right> = \left< D_X Y, Z \right> + \frac{1}{2} \left< (D_X
J) J Y, Z \right>.
\end{gather}

%%%%%%%%%%%%%%%%%%%%%%%%%%%%%%%%%%%%%%%%%%%%%%%%%%%%%%%%%%%%%%
\section{Curvature Identities for Almost Hermitian Structures}
\label{identities}
%%%%%%%%%%%%%%%%%%%%%%%%%%%%%%%%%%%%%%%%%%%%%%%%%%%%%%%%%%%%%%

In this section we collect some important identites for the curvature and
torsion of almost Hermitian pairs.  Fix $(\omega, J)$ an almost Hermitian pair,
and let $g$ denote the associated Riemannian metric.
As usual, let $\Ric$ denote the usual Ricci curvature of the Levi-Civita
connection, and let $\Ric^{J}$ denote the $J$-invariant part of the Ricci tensor
of
$g$, i.e.
\begin{align} \label{RicJsym}
\Ric^{J} = \frac{1}{2} \left[ \Ric(\cdot, \cdot) + \Ric (J \cdot, J \cdot)
\right].
\end{align}
Furthermore set
\begin{align} \label{rhodef}
\rho(\cdot, \cdot) =&\ \Ric^{J}(J \cdot, \cdot).
\end{align}
Note $\rho \in \Lambda^{1,1}$.  Next set
\begin{align} \label{rhostardef}
\rho^* = R(\omega)
\end{align}
i.e., the Levi-Civita curvature operator acting on the K\"ahler form $\omega$.
One can see \cite{Apostolov} for more information on these quantities.

Now let $\N$ denote a Hermitian connection associated to an almost Hermitian
pair.  The connection $\N$ induces a Hermitian connection on the
anticanonical bundle, and we denote the curvature form of this
connection by $P$.  Alternatively, if $\Omega$ denotes the curvature of $\N$,
one has
\begin{align} \label{Pdef}
P_{ij} = \omega^{kl} \Omega_{ij k l}.
\end{align}
By the general Chern-Weil theory, $P$ is a closed form and $P \in \pi c_1(M,
J)$.  We record some lemmas relating these
different curvature tensors.  A
key role is played in our analysis by the Weitzenb\"ock formula for two-forms.

\begin{lemma} \label{weitzenbock} Let $(M^{2n}, \omega, J)$ be an almost
Hermitian manifold.  Then
\begin{align} \label{curvaturecalc1}
{\rho^* - 2 \rho} =&\ {\left(D^* D \omega - \gD_d \omega \right).}
\end{align}
\begin{proof} By the Weitzenb\"ock formula for $2$-forms (\cite{Besse} pg. 53)
applied to $\omega$ we
conclude
\begin{align*}
\gD_d \omega - D^* D \omega =&\ \Ric(\omega \cdot, \cdot) - \Ric(\cdot, \omega
\cdot) -  {R(\omega)},
\end{align*}
where here the action of the Ricci tensor on the two form $\omega$ is by raising
the index on the Ricci tensor using the metric and letting the endomorphism act
naturally.  Phrasing this in terms of $J$ one sees
\begin{align*}
{R(\omega)} + \left[ \Ric(\cdot, J \cdot) - \Ric(J \cdot, \cdot) \right] = D^* D
\omega - \gD_d \omega.
\end{align*}
The Ricci curvature terms simplify to $- 2 \rho$, and the result follows.
\end{proof}
\end{lemma}

\noindent Furthermore (see \cite{Apostolov}), for an almost K\"ahler structure
one has the relation
\begin{align} \label{curvcalc2}
P =&\ \rho^* - \frac{1}{2} N^1,
\end{align}
where
\begin{align} \label{N1def}
N^1(X, Y) = \left< D_{JX} \omega, D_{Y} \omega \right>,
\end{align}
or, in coordinates,
\begin{align*}
N^1_{ab} =&\ g^{kl} g_{mn} J_a^p D_p J_k^m D_b J_l^n.
\end{align*}

\noindent We next derive a more general version of this formula.

\begin{lemma} \label{generalPformula} Let $(M^{2n}, \omega, J)$ be an almost
Hermitian manifold.  Let $\N$ denote the canonical connection corresponding to
$t = 0$ in the sense of (\ref{generalconnection}).
 Then
\begin{align*}
P =&\ \rho^* - \frac{1}{2} N^1 + \frac{1}{2} W
\end{align*}
where
\begin{align*}
W(X, Y) =&\ \left< [J D_X J, D_{JX} J ], D_Y J \right>.
\end{align*}
\begin{proof} Fix commuting vector fields $X, Y$, and let $e_i$ be local normal
coordinates for $g$ at some point.  Then
\begin{gather} \label{generalPloc10}
\begin{split}
P(X, Y) =&\ \Omega \left(X, Y, e_i, J e_i \right)\\
=&\ {\left< \N_Y \N_X e_i - \N_X \N_Y e_i, J e_i \right)}\\
=&\ {\left< D_Y \left( D_X e_i + \frac{1}{2} (D_X J) Je_i \right)+ \frac{1}{4}
(D_Y J) J(D_X J)(J e_i), J e_i \right>} \\
&\ - \mbox{symmetric term in $X$ and $Y$}\\
=&\ R(\omega) + {\frac{1}{2} \left< \left[ (D_Y D_X - D_X D_Y )J \right]J e_i, J
e_i \right>}\\
&\ {+ \frac{1}{2} \left< (D_X J) (D_Y J) e_i - (D_Y J) (D_X J) e_i, J
e_i\right>}\\
&\ {+
\frac{1}{4} \left< (D_Y J) J (D_X J) J e_i - (D_X J) J (D_Y J) J e_i, J e_i
\right>.}
\end{split}
\end{gather}

{First we observe that
\begin{align*}
g^{kl} \left< (D_i D_j J - D_j D_i J) J e_k, J e_l \right>
=&\ g^{kl} g_{st} \left(R_{ij r}^p J_p^s - R_{i j p}^s J_r^p \right) J_k^r
J_l^t\\
=&\ R_{ijr}^k J_k^r + g^{kl} g_{st} R_{ijk}^s J_l^t\\
=&\ R_{ijr}^k J_k^r - R_{ijk}^s J_s^k\\
=&\ 0.
\end{align*}
Therefore the second term in the last equality of (\ref{generalPloc10})
vanishes.  Next we note that
\begin{align*}
- \frac{1}{2} \left< (D_Y J) (D_X J) e_i , J e_i\right> =&\ \frac{1}{2} \left<
J(D_Y J)(D_X J) e_i, e_i \right>\\
=&\ - \frac{1}{2} \left< (D_Y J) J (D_X J) e_i, e_i \right>\\
=&\ + \frac{1}{2} \left< (D_Y J) \left( D_{JX} J + [J D_X J, D_{JX} J ] \right)
e_i, e_i \right>\\
=&\ - \frac{1}{2} \left< D_{JX} \omega, D_Y \omega \right> + \frac{1}{2} \left<
[J D_X J, D_{JX} J], D_Y J \right>.
\end{align*}
A similar calculation yields the same result for the skew symmetric term.  For
the remaining term we compute
\begin{align*}
- \frac{1}{4} \left< (D_X J) J (D_Y J) J e_i, J e_i \right> =&\ \frac{1}{4}
\left< J (D_X J) (D_Y J) J e_i, J e_i \right>\\
=&\ - \frac{1}{4} \left< \left( D_{JX} J + [J D_X J, D_{JX} J] \right)  (D_Y J)
Je_i, Je_i \right>\\
=&\ \frac{1}{4} \left< D_{JX} \omega, D_Y \omega \right> - \frac{1}{4} \left<
[J D_X J, D_{JX} J ], D_Y J \right>
\end{align*}
A similar calculation yields the same result for the skew symmetric piece, and
the result follows.}
\end{proof}
\end{lemma}

It is relevant to us to know that the commutator term in the definition of $W$
above is determined by $d \omega$.  We record the formula here.

\begin{lemma} \label{generalJlemma} Let $(M^{2n}, \omega, J)$ be an almost
Hermitian manifold.  Then
\begin{align} \label{generalJformula}
D_{JX} \omega(Y, Z) - D_X \omega(JY, Z) = \left( d \omega \right)^+(J X, Y, Z) -
\left(d \omega \right)^+(JX, JY, JZ).
\end{align}
In particular, if $d \omega = 0$ then one has
\begin{align} \label{AKJformula}
D_{JX} J = - J(D_X J).
\end{align}
\begin{proof} This is a restatement of \cite{GauduchonConn} Proposition 1.iv.
\end{proof}
\end{lemma}

\begin{lemma} \label{P2formula} Let $(M^{2n}, \omega, J)$ be an almost
Hermitian manifold.  Let $\N$ denote the canonical connection corresponding to
$t = 0$ in the sense of (\ref{generalconnection}).  Then
\begin{align} \label{P20part}
P^{(2,0) + (0,2)} =&\ {D^* D \omega - N^2 + \left( - \gD_d \omega
- \frac{1}{2} N^1 + \frac{1}{2} W \right)^{(2,0) + (0,2)}}
\end{align}
where
\begin{align} \label{N2def}
N^2(X, Y) =&\ \left< (D J) JX, (DJ) Y \right>,
\end{align}
or, in coordinates,
\begin{align*}
N^2_{a b} =&\ g^{ij} g_{mn} D_i J_p^m J_a^p D_j J_b^n.
\end{align*}
\begin{proof} Combining Lemmas \ref{weitzenbock} and \ref{generalPformula}
yields
\begin{align} \label{Pformula}
P = {2 \rho + D^* D \omega - \gD_d \omega - \frac{1}{2} N^1 + \frac{1}{2} W}
\end{align}
for an almost K\"ahler structure.  Since $\rho \in \Lambda^{1,1}$, it remains to
compute the $(2, 0) + (0,2)$
component of $D^* D \omega$.  We do this by computing the $(1,1)$ component,
which we will compute in local coordinates.
\begin{align*}
- \left( D^* D \omega \right)^{1,1}_{ab} =&\ - \frac{1}{2} \left[ \left( D^* D
\omega \right)(J, J) + D^* D \omega \right]_{ab}\\
=&\ \frac{1}{2} g^{ij} \left[ (D_i D_j \omega_{pq}) J^p_a J^q_b + D_i D_j
\omega_{ab} \right]\\
=&\ \frac{1}{2} g^{ij} \left[ D_i D_j \left( \omega_{pq} J^p_a J^q_b \right) +
D_i D_j \omega_{ab} {- 2 D_i \omega_{pq} D_j J_a^p J_b^q - 2 D_i \omega_{pq}
J_a^p
D_j J_b^q} \right.\\
&\ \qquad \left. - \omega_{pq} \left( (D_i D_j J^p_a) J^q_b + D_i J^p_a D_j
J^q_b + D_j J^p_a D_i
J^q_b + J^p_a D_i D_j J^q_b \right) \right].
\end{align*}
Using compatibility of $\omega$ with $J$,
\begin{align*}
D_i D_j \left( \omega_{pq} J^p_a J^q_b \right) = D_i D_j \omega_{ab}.
\end{align*}
Also, we have that
\begin{align*}
- \omega_{pq} (D_i D_j J^p_a) J^q_b =&\ - g_{pb} D_i D_j J^p_a\\
=&\ - D_i D_j \left( g_{pb} J^p_a \right)\\
=&\ - D_i D_j \left( \omega_{ab} \right).
\end{align*}
Next we compute
\begin{align*}
 - \omega_{pq} J^p_a D_i D_j J^q_b =&\ g_{aq} D_i D_j J^q_b\\
=&\ - D_i D_j \omega_{ab}.
\end{align*}
{Next note that
\begin{align*}
D_i \omega_{pq} D_j J_a^p J_b^q =&\ D_i \left[ - J_q^r g_{pr} \right] D_j J_a^p
J_b^q\\
=&\ - g_{pr} \left[ D_i J_q^r J_b^q \right] D_j J_a^p\\
=&\ g_{pr} J_q^r D_i J_b^q D_j J_a^p\\
=&\ - \omega_{pq} D_j J_a^p D_i J_b^q.
\end{align*}
Likewise $D_i \omega_{pq} J_a^p D_j J_b^q = - \omega_{pq} D_i J_a^p D_j J_b^q$.}
It follows that
\begin{align*}
 \left( D^* D \omega \right)^{1,1}_{ab} = {- g^{ij} \omega_{pq} D_i J^p_a D_j
J^q_b} = { g^{ij} g_{mn} D_i J_p^m J_a^p D_j J_b^n.}
\end{align*}
The lemma follows.
\end{proof}
\end{lemma}

%%%%%%%%%%%%%%%%%%%%%%%%%%%%%%%%%%%%%%%%%%%%%%%%%%%
\section{Variations of Almost Hermitian Structures} \label{variations}
%%%%%%%%%%%%%%%%%%%%%%%%%%%%%%%%%%%%%%%%%%%%%%%%%%%

\begin{lemma} \label{Kcondition} Let $(M^{2n}, J)$ be a complex manifold and
suppose $J(t)$ is a one-parameter family of endomorphisms of $TM$ such that
$J(0) = J$.  Then $J(t)$ is a one-parameter family of almost-complex structures
if and only if for all $t$,
\begin{align*}
J \left( \dt J \right) + \left( \dt J \right) J = 0.
\end{align*}
\begin{proof} Assume $J(t)$ is an almost complex structure.  Then
$J(t)^2 = - \Id$, thus
\begin{align*}
0 =&\ \dt J^2 = J K + K J.
\end{align*}
Conversely, if this equation holds for all time, we may integrate it to obtain
that $J^2(t) = J^2(0) = - \Id$, and so $J(t)$ is a one parameter family of
almost complex structures.
\end{proof}
\end{lemma}

\begin{lemma} \label{compatibilityconds} Let $(M^{2n}, g(t), J(t))$ be a
one-parameter family of metrics
and almost complex structures, with $g(0)$ compatible with $J(0)$.  Further
suppose
\begin{align*}
\dt g =&\ H + F\\
\dt J =&\ K
\end{align*}
where $F \in \Sym^{(2,0) + (0,2)} T^* M$ and $H \in \Sym^{1,1} T^* M$.  Then
$g(t)$ is compatible with $J(t)$ if and only if
\begin{align*}
F =&\ \frac{1}{2} \left( g(K, J) + g(J, K) \right).
\end{align*}
Furthermore, assuming this holds one has
\begin{align*}
\dt \omega = H(J \cdot, \cdot) + \frac{1}{2} \left[ g(K \cdot, \cdot) - g(\cdot,
K \cdot) \right]
\end{align*}
\begin{proof} It suffices to show that the time derivative of the compatibility
condition vanishes at time $t = 0$.  We directly compute
\begin{align*}
\dt \left( g(J \cdot, J \cdot) - g( \cdot, \cdot) \right) =&\ \left(H +
F \right)(J \cdot, J \cdot) + g(K \cdot, J \cdot) + g(J \cdot, K \cdot)\\
&\ - \left( H + F \right) (\cdot, \cdot).
\end{align*}
Note that $H(J \cdot, J \cdot) - H(\cdot, \cdot) = 0$.  Now let $F =
\frac{1}{2} \left( g(K \cdot, J \cdot) + g(J \cdot, K \cdot) \right)$.  Observe
that
\begin{align*}
F(J \cdot J \cdot) - F(\cdot, \cdot) =&\ \frac{1}{2} \left[ g(K J, JJ) + g(JJ,
KJ) - g(K, J) - g(J, K) \right]\\
=&\ \frac{1}{2} \left[ - g(K J, \cdot) - g(\cdot, KJ) - g(K, J) - g(J, K)
\right]
\end{align*}
Using Lemma \ref{Kcondition} and compatibility of $J$ and $g$ at time zero we
note that
\begin{align*}
- g(KJ, \cdot) = g(J K, \cdot) = g(JJ K, J) = - g(K, J)
\end{align*}
and likewise $- g(\cdot, KJ) = - g(J, K)$.  Thus combining these calculations we
conclude that the time derivative vanishes, and the lemma follows.
\end{proof}
\end{lemma}

\begin{lemma} \label{formcompatibility} Let $(M^{2n}, \omega(t), J(t))$ be a
one-parameter family of
K\"ahler forms and almost-complex structures with $\omega(0)$ compatible with
$J(0)$.  Further suppose
\begin{align*}
\dt \omega =&\ \phi + \psi\\
\dt J =&\ K
\end{align*}
where $\psi \in \Lambda^{(2,0) + (0,2)}$ and $\phi \in \Lambda^{1,1}$.  Then
$\omega(t)$
is compatible with $J(t)$ if and only if
\begin{align} \label{compatibilitycond}
\psi = \frac{1}{2} \left[ \omega(K, J) + \omega(J, K) \right].
\end{align}
\begin{proof} We directly compute
\begin{align*}
0 =&\ \dt \left( \omega(J \cdot, J \cdot) - \omega(\cdot, \cdot) \right)\\
=&\ \left( \phi + \psi\right) (J \cdot, J \cdot) + \omega(K \cdot, J \cdot) +
\omega(J \cdot, K \cdot) - \left(\phi + \psi
\right) (\cdot, \cdot)\\
=&\ \psi (J \cdot, J \cdot) - \psi(\cdot, \cdot) +
\omega(K \cdot, J \cdot) + \omega(J \cdot, K \cdot).
\end{align*}
Since $\psi \in \Lambda^{(2,0) + (0,2)}$ the above equation is equivalent to
\begin{align*}
2 \psi = \omega(K, J) + \omega(J, K)
\end{align*}
as required.
\end{proof}
\end{lemma}

\begin{rmk} \label{comprmk}
Fix a point $p \in M$ and choose some local coordinates.  Certainly
(\ref{compatibilitycond}) holds if
\begin{align*}
K_a^b = g^{b c} \psi_{a c}.
\end{align*}
Observe however that $K$ is not
determined by $\psi$ alone.  Indeed one may add to $K$ an endomorphism of the
form $g^{-1} W^{-J}$, where $W$ is a symmetric two tensor,
and (\ref{compatibilitycond}) will still hold.
\end{rmk}

\begin{lemma} \label{liederivative} Let $(M^{2n}, J)$ be an almost-complex
manifold and let $X$ be a
vector field on $M$.  Then
\begin{align} \label{Jlieskew}
J L_X J + L_X J J = 0.
\end{align}
Furthermore, if $\omega$ is compatible with $J$, we have
\begin{align*}
\left( L_X \omega \right)^{(2,0) + (0,2)} =&\ \frac{1}{2} \left( \omega( L_X J
\cdot, J
\cdot) + \omega(
J \cdot, L_X J \cdot ) \right).
\end{align*}
\begin{proof} By \cite{Besse} pg. 86, we have the formula for $L_X J$:
\begin{align} \label{lieJ}
\left(L_X J \right)(Y) =&\ [X, J Y] - J[X, Y].
\end{align}
Given this,
the first equation follows by direct calculation.  The second equation obviously
must hold since it is just the linearized compatibility condition
(\ref{compatibilitycond}) and the action of a diffeomorphism preserves
compatibility, but we just as well compute
\begin{align*}
0 =&\ L_X \left( \gw(\cdot, \cdot) - \gw( J \cdot, J \cdot) \right)\\
=&\ (L_X \gw)(\cdot, \cdot) - \left(L_X \gw \right)(J \cdot, J \cdot) - \gw( L_X
J
\cdot, J \cdot) - \gw( J \cdot, L_X J \cdot).
\end{align*}
Rearranging the above formula gives the result.
\end{proof}
\end{lemma}

\begin{lemma} \label{coordianteLieJ} Let $(M^{2n}, J)$ be an almost complex
manifold and let $X$ be a vector field on $M$.  Then
\begin{align*}
\left( L_X J\right)_k^l = J_p^l \del_k X^p - J_k^p \del_p X^l + X^p
\del_p J_k^l.
\end{align*}
\begin{proof} Choose local coordinate vector fields $e^k$.  Using (\ref{lieJ})
we see
\begin{align*}
\left( L_X J \right)_k^l e^k =&\ - \left( J e^k \right)^p \del_p X^l +
J_p^l \left[ e^k  \del_k X^p \right] + X^p \del_p J_k^l\\
=&\ - J_k^p \del_p X^l + J_p^l \del_k X^p + X^p \del_p J_k^l.
\end{align*}
as required.
\end{proof}
\end{lemma}

%%%%%%%%%%%%%%%%%%%%%%%%%%%%%%%%%%%%%%%%%%%%%%%%%%%%%%%%
\section{Parabolic Flows of Almost Hermitian structures} \label{generalste}
%%%%%%%%%%%%%%%%%%%%%%%%%%%%%%%%%%%%%%%%%%%%%%%%%%%%%%%%

In this section we prove Theorem \ref{almosthermitianste}.  Let us recall some
definitions from the introduction used in (\ref{almostcomplexflow}).  In
particular, let $(M^{2n}, \gw, J)$ be an almost Hermitian manifold and let $\N$
denote the associated Chern connection (see Definition \ref{Cherndef}).  Let
$\Omega$ denote the $(4,0)$ curvature tensor associated to $\N$, and consider
\begin{align*}
S_{ij} = \omega^{kl} \Omega_{klij}.
\end{align*}
Furthermore, let $N$ denote the Nijenhuis tensor associated to $J$, and let
\begin{gather*}
\mathcal K_j^i = \omega^{kl} \N_k N_{l j}^i.
\end{gather*}
\noindent Let $Q$ denote a $(1,1)$ tensor which is quadratic expression in the
torsion of $\N$, and let $\mathcal H$ denote a $J$-skew endomorphism of the
tangent bundle which again is quadratic in the Nijenhuis tensor.  Let
\begin{align} \label{Fdef}
H =&\ \frac{1}{2} \left[ \gw(\mathcal K - \mathcal H, J) + \gw(J, \mathcal K -
\mathcal H) \right].
\end{align}
Consider the initial value problem
\begin{gather} \label{flow}
\begin{split}
\dt \omega =&\ - S + Q + H\\
\dt J =&\ - \mathcal K + \mathcal H\\
\omega(0) =&\ \omega_0\\
J(0) =&\ J_0.
\end{split}
\end{gather}

\noindent We first show some preliminary lemmas which show that the right hand
sides of (\ref{flow}) satisfy the linearized conditions for one parameter
families of almost Hermitian pairs from Lemmas \ref{Kcondition} and
\ref{formcompatibility}

\begin{lemma} \label{Nijenhuissymmetry} Let $(M^{2n}, J)$ be an almost-Hermitian
manifold.  Then, viewing the Nijenhuis tensor $N$ as a section of $\Lambda^1
\otimes \End(TM)$,
\begin{align*}
J N + N J = 0.
\end{align*}
\begin{proof} We can derive this by direct calculation using the definition of
the Nijenhuis tensor (\ref{Nijenhuis}).  First
\begin{align*}
-J N =&\ J \left( [X, Y] + J \left( [J X, Y] + [X, J Y] \right) - [JX, JY]
\right)\\
=&\ J[X, Y] - [JX, Y] - [X, J Y] - J [JX, JY].
\end{align*}
Next
\begin{align*}
-N J =&\ [X, JY] + J \left( [JX, JY] + [X, JJY] \right) - [JX, JJY]\\
=&\ [X, JY] + J[JX,JY] - J[X, Y] + [JX, Y].
\end{align*}
The result follows adding these two calculations together.
\end{proof}
\end{lemma}

\begin{lemma} \label{Pcompatibility} Let $(M^{2n}, \gw, J)$ be an almost
Hermitian
manifold.  Then
\begin{align*}
J  \PP + \PP J = 0.
\end{align*}
\begin{proof} We may write the result of Lemma \ref{Nijenhuissymmetry} in
coordinates as
\begin{align*}
J_k^m N_{j m}^l + N_{j k}^m J_m^l = 0
\end{align*}
We differentiate this using the Chern connection.  Since $J$ is parallel we see
\begin{align*}
0 =&\ J_k^m \N_i N_{jm}^l + \N_i N_{jk}^m J_m^l
\end{align*}
we can now take the required contraction of indices using $\omega$ to yield the
statement of the lemma.
\end{proof}
\end{lemma}

\begin{defn} \label{Xdefn} Let $(M^{2n}, \gw, J)$ be an almost Hermitian
manifold.
 Let
$\bar{\N}$ denote some fixed connection on $TM$.  Define a vector field
\begin{gather} \label{Xdef}
X^p = X(\omega, J, \bar{\N})^p = \omega^{k l} \bar{\N}_k J_{l}^p.
\end{gather}
\end{defn}

\begin{prop} \label{Joperator} Let $(M^{2n}, \gw, J)$ be an almost-Hermitian
manifold and let
$\bar{\N}$ denote some fixed connection on $TM$.  The map
\begin{align*}
J \to \mathcal K - L_{X(\omega, J, \bar{\N})} J
\end{align*}
is a second order elliptic operator.
\begin{proof}  We recall a coordinate formula for the Nijenhuis tensor.
\begin{align} \label{Nijencoords}
N_{j k}^{i} =&\ J_j^p \del_p J_k^i - J_k^p \del_p J_j^i - J_p^i \del_j J_k^p +
J_p^i \del_k J_j^p.
\end{align}
It follows that
\begin{gather} \label{Kcalc}
\begin{split}
\mathcal K_j^i =&\ \omega^{kl} \N_k N_{l j}^i\\
=&\ \omega^{kl} \left( J_l^q \del_k \del_q J_j^i - J_j^q \del_k \del_q J_l^i -
J_q^i \del_k \del_l J_j^q + J_q^i \del_k \del_j J_l^q \right) + \mathcal O(\del
J, \del \omega)\\
=&\ - g^{q k} \del_k \del_q J_j^i - \omega^{kl} \left( J_j^q \del_k \del_q J_l^i
+ J_q^i \del_k \del_l J_j^q - J_q^i \del_k \del_j J_l^q \right) + \mathcal
O(\del J, \del \omega)
\end{split}
\end{gather}
where the notation $\mathcal O(\del J, \del \omega)$ means an expression which
only depends
on at most first derivatives of $J$ and $\omega$ (possibly in a nonlinear
fashion).  In particular, Chern connection terms are of this form.  Note that
the matrix $\gw$ is skew-symmetric, but
coordinate derivatives are symmetric, therefore the middle term in the
parentheses in the last line vanishes.  Also, using (\ref{Xdef}) and Lemma
\ref{coordianteLieJ} we express
\begin{align*}
\left[L_{X(\gw, J, \bar{\N})} J \right]_j^i =&\ \omega^{kl} \left( J_q^i \del_j
\del_k J_l^q  - J_j^q \del_q
\del_k J_l^i \right) + \mathcal O(\del J, \del
\omega).
\end{align*}
Combining these two calculations yields
\begin{gather*}
\begin{split}
\left[ \mathcal K - L_{X(\gw, J, \bar{\N})} J \right]_j^i = - g^{k l} \del_k
\del_l
J_j^i + \mathcal O(\del J, \del \gw).
\end{split}
\end{gather*}
The claim follows immediately.
\end{proof}
\end{prop}

\begin{prop} \label{goperator} Let $(M^{2n}, \gw, J)$ be an almost-Hermitian
manifold.  The map
\begin{align*}
\omega \to  S(\omega)
\end{align*}
is a second order elliptic operator.
\begin{proof} Fix a point $p \in M$, and choose a local frame of $(1,0)$ vector
fields $\{e_i \}$ such that $g_{i \bj}(p) = \gd_{ij}$.  In this frame we compute
using metric compatibility of $\N$,
\begin{align*}
S_{k \bl} =&\ \gw^{i \bj} \Omega_{i \bj k \bl}\\
=&\ \gw^{i \bj} \left< \N_{i} \N_{\bj} e_k - \N_{\bj} \N_i e_k - \N_{[e_i,
e_{\bj}]} e_k, e_{\bl} \right>\\
=&\ \gw^{i \bj} \left( e_i \left< \N_{\bj} e_k, e_{\bl} \right> - \left<
\N_{\bj}
e_k, \N_i e_{\bl} \right> - e_{\bj} \left< \N_i e_k, e_{\bl} \right> + \left<
\N_i e_k, \N_{\bj} e_{\bl} \right> \right) + \mathcal O(\del \gw, \del J)\\
=&\ \gw^{i \bj} \left( e_i \left< \N_{\bj} e_k, e_{\bl} \right> - e_{\bj} e_{i}
\left< e_{k}, e_{\bl} \right> + e_{\bj} \left< e_k, \N_i e_{\bl} \right>
\right) + \mathcal O(\del \gw, \del J).
\end{align*}
Now using $J$ compatibility of the connection and the fact that the torsion $T$
has no $(1,1)$-component we see that
\begin{align*}
\left< \N_{\bj} e_k, e_{\bl} \right> =&\ \left< \N_{\bj} e_k - \N_{k} e_{\bj},
e_{\bl} \right>\\
=&\ \left< T_{\bj k} + [e_{\bj}, e_k], e_{\bl} \right>\\
=&\ \left<[e_{\bar{j}}, e_k], e_{\bl} \right>\\
=&\ \mathcal O(\gw, \del J).
\end{align*}
The last line follows since the basis $e_i$ is constructed by projecting
local coordinates onto $T^{1,0}$ using $J$, and therefore their Lie brackets
will only contain derivatives of $J$.  Likewise one concludes that $\left< e_k,
\N_{e_i} e_{\bl}
\right> = \mathcal
O(\gw, \del J)$.  Therefore
\begin{align*}
S_{k \bl} =&\ - g^{i \bj} e_{\bj} e_i \gw_{k \bl} + \mathcal O(\del \gw, \del^2
J)\\
=&\ - \frac{1}{2} g^{ab} \del_a \del_b \gw_{k \bl} + \mathcal O(\del \gw, \del^2
J).
\end{align*}
The result follows.
\end{proof}
\end{prop}

\noindent We can now give the proof of Theorem \ref{almosthermitianste}.
\begin{proof} First we show existence.  This will be a two step process.  First
we will define a gauge-fixed flow which will define a strictly parabolic system.
 This flow equation will only be defined however for compatible pairs.
Therefore to apply short time existence results from the theory of parabolic
differential equations one needs to define a more general evolution equation
which makes sense for arbitrary pairs $(\omega, J)$.  Moreover, the desired flow
on $J$ takes place in a nonlinear manifold, therefore one must ``pull back'' the
flow on $J$ to a linear space, namely the tangent space to the space of almost
complex structures at $J_0$.  We define such a generalized version of
our gauge-fixed flow which has short time existence.  We can show that this flow
preserves compatibility of the initial condition and produces a solution of the
original gauge-fixed flow.  We remove the gauge parameter to finally produce the
required solution of the original flow.

Fix $\bar{\N}$ any connection on
$TM$ and let $X$ be defined as in Definition \ref{Xdef}.  First consider the
following
gauge-fixed version of equation (\ref{flow})
\begin{gather} \label{gfflow}
\begin{split}
\dt \omega =&\ - S + Q + H + L_{X(g,J)} \gw = \mathcal D_1(\gw,J)\\
\dt J =&\ - \mathcal K + \mathcal H + L_{X(g,J)} J = \mathcal D_2(\gw,J).
\end{split}
\end{gather}
We observe that by definition the vector field $X(g,J)$ can be expressed
completely in terms
of first derivatives of $J$ and therefore $\left(L_{X(\gw,J)} \gw\right)_{ij}$
is a first order operator in $\omega$.  Use $\mathcal L_{\gw}, \mathcal L_{J}$
to
denote linearization in the $\gw$ and $J$ variables respectively.  It follows
from Proposition
\ref{goperator} that
\begin{align*}
\gs \left[ \widehat{\mathcal L_{\gw} \mathcal D_1} \right](h)_{ij} =&\ \gs
\left[
\widehat{\mathcal L
(-2S)} \right](h)_{ij}\\
=&\ \brs{\xi}^2 h_{i j}.
\end{align*}
Furthermore, from Proposition \ref{Joperator} we conclude that
\begin{align*}
\gs \left[ \widehat{\mathcal L_{J} \mathcal D_2} \right](K)_i^j =&\ \brs{\xi}^2
K_i^j
\end{align*}
We also need to check the linearization of $\mathcal D_2$ in the variable
$\gw$.
Since by construction we have that $\mathcal D_2$ only depends on first
derivatives of $\gw$, we conclude
\begin{align*}
\gs \left[ \widehat{\mathcal L_{\gw} \mathcal D_2} \right](h)_i^j =&\ 0.
\end{align*}
We note that second derivative terms of $J$ appear in the evolution of $\gw$,
therefore these terms appear in the full linearized operator.  Collecting these
observations we conclude that the overall symbol is upper-triangular.  In
particular it takes the form
\begin{align*}
\gs \left[ \widehat{\mathcal L \mathcal D} \right](h, K) =&\
\left(
\begin{matrix}
I & *\\
0 & I
\end{matrix}
\right) \left(
\begin{matrix}
h\\
K
\end{matrix}
\right)
\end{align*}
It follows that (\ref{gfflow}) is a strictly parabolic system of equations.
However, as mentioned above we must define a more general flow defined for
arbitrary pairs $(\omega, J)$ to apply short time existence theory.
First note that the
space of almost complex structures $\mathbb J$ near a fixed $J_0$ is a smooth
Banach manifold with tangent space modeled on the space of $J_0$-skew
endomorphisms (see Lemma \ref{Kcondition}).  On obtains a diffeomorphism
\begin{align*}
\pi : U_0 \subset \mathbb T \mathbb J_{J_0} \to U_{J_0} \subset \mathbb J
\end{align*}
from a neighborhood of $0$ in $\mathbb T \mathbb J_{J_0}$ to a neighborhood of
$J_0$ in $\mathbb J$.  Moreover this diffeomorphism satisfies
\begin{align} \label{Dpicalc}
D \pi_0 = \Id_{\mathbb T \mathbb J_{J_0}}.
\end{align}
Note that our desired flow for $J$ is moving through a certain (nonlinear)
manifold in the space of endomorphisms of the tangent bundle.  We will use the
map $\pi$ to pull back the flow onto the linear space $\mathbb T \mathbb J$.

Suppose now $(\omega_0, J_0)$ is a compatible
pair and let $g_0$ denote the associated metric.  Recall that if $J$ is an
almost complex structure and $\omega \in \Lambda^2$, $\omega^J$
denotes the $J$-symmetric piece of $\omega$, while $\omega^{-J}$ denotes the
$J$-anti-invariant piece.  Consider now the initial value
problem
\begin{gather} \label{gfmodflow}
\begin{split}
\dt \omega =&\ \mathcal D_1(\gw^{\pi E}, \pi E) - D^*_{g_0} D_{g_0} \left(
\omega^{- \pi E} \right) =: \widetilde{\mathcal D}_1(\omega, E)\\
\dt E =&\ \left(D \pi_{\pi E}^{-1}\right) \left(\mathcal D_2(\omega^{\pi E}, \pi
E) \right) =:\widetilde{\mathcal D}_2(\omega, E).\\
\omega(0) =&\ \omega_0\\
E(0) =&\ 0.
\end{split}
\end{gather}
We observed above that $\mathcal D_2(\omega^{\pi E}, \pi E)$ lies in $\mathbb T
\mathbb J_{\pi E}$, therefore the operator $\widetilde{\mathcal D_2}$ is well
defined, and has image in $\mathbb T \mathbb ( \mathbb T \mathbb J)_{E} \cong
\mathbb T \mathbb J_{J_0}$ for arbitrary pairs
$(\omega, J)$, therefore this equation defines a flow in the linear space $B :=
\mathbb T \mathbb J_{J_0} \oplus \Lambda^2_{\mathbb R}$.

We want to compute the linearization of this system at $t = 0$.  First we
compute the linearization of $\widetilde{\mathcal D}_1$ in the $\omega$
variable.  Combining Proposition \ref{goperator} with
an obvious calculation of the symbol for $D^*_{g_0} D_{g_0} ( \omega^{- \pi E})$
yields
\begin{align*}
\sigma \left[ \mathcal L_{\omega} \left( \mathcal D_1(\omega^{\pi E}, \pi E) -
D^*_{g_0} D_{g_0} \left( \omega^{- \pi E} \right) \right)
\right]^{\wedge}(h)_{ij} =
\brs{\xi}^2 h_{ij}^{\pi E} +  \brs{\xi}^2 h^{- \pi E}_{ij} = \brs{\xi}^2 h_{ij}.
\end{align*}
Furthermore, the calculation of the linearization of $\widetilde{\mathcal
D}_2(\omega, E)$ at $t = 0$ is identical to that for $\mathcal D_2$
above using (\ref{Dpicalc}).  It follows that
\begin{align*}
\gs \left[ \widehat{\mathcal L \widetilde{\mathcal D}} \right](h, K) =&\
\left(
\begin{matrix}
I & *\\
0 & I
\end{matrix}
\right) \left(
\begin{matrix}
h\\
K
\end{matrix}
\right).
\end{align*}
Therefore the initial value problem (\ref{gfmodflow}) is a nonlinear strictly
parabolic equation in the Banach space $B$, and standard
results imply the existence of a short time solution to (\ref{gfmodflow}).

Now let $J = \pi E$.  We claim that $(\omega, J)$ this is in fact a
solution to (\ref{gfflow}).  First we compute
\begin{align*}
\dt J = \left(D \pi_{E} \right) \left(\dt J \right) =\left(D \pi_{E} \right)
\left( D \pi_{\pi E}^{-1} \right) \left(\mathcal D_2(\omega^{J}, J)
\right) = \mathcal D_2(\omega^{J}, J).
\end{align*}
Next we want to show
compatibility of the pair $(\omega, J)$ is preserved.  Note that since we have
already computed that solutions to (\ref{gfflow}) satisfy the conditions of
Lemma \ref{formcompatibility} one has
\begin{align*}
\dt \omega^{- J} =&\ \left[- D^*_{g_0} D_{g_0} \left( \omega^{- J} \right)
\right]^{-J}.
\end{align*}
Now let $\til{g}(t)$ be a one-parameter family of metrics which is compatible
with $J(t)$.
It follows that
\begin{gather} \label{genpfloc10}
\begin{split}
\dt \brs{\omega^{-J}}_{\til{g}} =&\ 2 \left< - D^*_{g_0} D_{g_0} \omega^{-J},
\omega^{-J} \right>_{g_0} + \left(
\dt \til{g} \right) * \left(\omega^{-J} \right)^{*2}\\
=&\ 2 \left< \tr_{g_0} \left( D^2_{\til{g}} + \til{R} + R_0 \right) \omega^{-J},
\omega^{-J} \right> + \left( \dt \til{g} \right) * \left(\omega^{-J}
\right)^{*2}\\
\leq&\ \tr_{g_0} D^2_{\til{g}} \brs{\omega^{-J}}_{g_0} - 2 \brs{D_{\til{g}}
\omega^{-J}}_{g_0, \til{g}}^2 + C \brs{\omega^{-J}}^2_{\til{g}}.
\end{split}
\end{gather}
Applying the maximum principle to $e^{-C t} \brs{\omega^{-J}}_{\til{g}}^2$ we
conclude that if $\omega^{-J}(0) \equiv 0$, then
$\omega^{-J}(t) \equiv 0$ for all $t$.  Hence the pair $(\omega(t), J(t))$ is
compatible for all $t$, and we conclude that the one parameter family
$(\omega(t), J(t))$ is a
solution to (\ref{gfflow}).

Now we want to
pull back our solution to (\ref{gfflow}) by the family of diffeomorphisms
generated by $X$.
Specifically let $\phi_t$ be a one-parameter family of diffeomorphisms of $M$
defined by the ODE
\begin{gather} \label{Xode}
\begin{split}
\dt \phi_t =&\ - X(\omega(t), J(t), \bar{\N})\\
\phi_0 =&\ \id_M.
\end{split}
\end{gather}
It follows that
\begin{align*}
\dt \left( \phi_t^* \gw(t) \right) =&\ \frac{\del}{\del s}|_{s = 0} \left(
\phi_{t+s}^* \gw(t+s) \right)\\
=&\ \phi_t^* \left( \dt \gw(t) \right) + \frac{\del}{\del s}|_{s = 0} \left(
\phi_{t+s}^* \gw(t) \right)\\
=&\ \phi_t^* \left( - S + Q + H + L_{X(\gw(t), J(t))} \gw \right) +
\frac{\del}{\del s} |_{s = 0} \left[ \left(\phi_t^{-1} \circ \phi_{t+s}
\right)^* \phi_t^* \gw_t \right]\\
=&\ \left( - S + Q + H \right) \left( \phi_t^*(\gw), \phi_t^*(J)
\right) + \phi_t^* \left( L_{X(\gw(t), J(t))} \right) - L_{ \left(\phi_t^{-1}
\right)_*
X(\gw(t), J(t))} \left( \phi_t^* \gw(t) \right)\\
=&\ \left( - S + Q + H \right) \left( \phi_t^*(\gw), \phi_t^*(J)
\right).
\end{align*}
Likewise we may compute
\begin{align*}
\dt \left( \phi_t^* J(t) \right) =&\ \frac{\del}{\del s} |_{s = 0} \left(
\phi_{t+s}^* J(t+s) \right)\\
=&\ \phi_t^* \left( \dt J(t) \right) + \frac{\del}{\del s}|_{s = 0} \left(
\phi_{t+s}^* J(t) \right)\\
=&\ \phi_t^* \left( - \PP(\gw, J) + \mathcal H(\gw, J) + L_{X(\gw(t), J(t))}
\right) +
\frac{\del}{\del s} |_{s = 0} \left[ \left(\phi_t^{-1} \circ \phi_{t+s}
\right)^* \phi_t^* J_t \right]\\
=&\ - \mathcal K(\phi_t^* \gw(t), \phi_t^* J(t)) + \mathcal H(\phi_t^* \gw(t),
\phi_t^* J(t))\\
&\ \qquad + \phi_t^* L_{X(\gw(t),
J(t))} - L_{(\phi_t^{-1})_* X(\gw(t), J(t))} \left(\phi_t^* J(t) \right)\\
=&\ - \mathcal K(\phi_t^* \gw(t), \phi_t^* J(t)) + \mathcal H(\phi_t^* \gw(t),
\phi_t^* J(t)).
\end{align*}
Therefore $\left(\phi_t^* \gw(t), \phi_t^*(J(t)) \right)$ is a solution to
(\ref{flow}).

Next we show uniqueness.  As in the proof of uniqueness for Ricci flow, we will
show that the diffeomorphism ODE (\ref{Xode}), when written with respect to
the changing metric, is in fact a parabolic equation for $\phi$.  What is more,
as we now show, our choice of $X$ is essentially equivalent to that used for
Ricci flow short-time existence.  Let $\gG_C, \gG, \bar{\gG}$ denote the
connection coefficients of the Chern , Levi-Civita, and background connections
respectively.  Consider the following calculation:
\begin{gather} \label{Xcalc}
\begin{split}
X^p =&\ \gw^{kl} \bar{\N}_k J_l^p\\
=&\ \gw^{kl} \del_k J_l^p + \mathcal O(\gw, J)\\
=&\ \gw^{kl} \left( \N_k J_l^p + \left(\gG_C \right)_{k l}^q J_q^p - \left(\gG_C
\right)_{k q}^p J_l^q
\right) + \mathcal O(\gw, J)\\
=&\ \gw^{kl} \left(\gG_C \right)_{kl}^q J_q^p + g^{kq} \left(\gG_C
\right)_{kq}^p + \mathcal O(\gw, J).
\end{split}
\end{gather}
The first term is the contraction of the Chern connection coefficient with a
skew-symmetric one-form, and hence vanishes.  Specifically we compute
\begin{gather} \label{Xcalc2}
\begin{split}
\omega^{kl} \left(\gG_C \right)_{kl}^q J_q^p
=&\ \frac{1}{2} \omega^{k l} \left( \left(\gG_C \right)_{kl}^q - \left(\gG_C
\right)_{lk}^q \right) J_q^p\\
=&\ \frac{1}{2} \omega^{k l} T_{k l}^q J_q^p\\
=&\ 0
\end{split}
\end{gather}
since the torsion of $\N$ has no $(1,1)$ component.  Next observe that
since $g$ is symmetric, the contraction $g^{k q} \gG_{k q}^p$
does note involve the torsion of the connection.  In particular we conclude
\begin{align*}
g^{kq} \left(\gG_C \right)_{kq}^p =g^{kq} \gG_{kq}^p
\end{align*}
In particular, combining these calculations we may
conclude that
\begin{align} \label{Xform}
X^p =&\ g^{k l} \left[ \gG_{k l}^p - \bar{\gG}_{kl}^p \right] + \mathcal O(\gw,
J).
\end{align}
In particular we have shown that, up to lower order terms, the vector field we
used in our short-time existence proof is the same as that used for Ricci flow.
So, set $\til{g} = \phi_t^* g(t)$, $\til{J} = \phi_t^* J(t)$.  It follows (see
\cite{Chow} pg. 89) that one may rewrite the solution to (\ref{Xode}) as
\begin{gather} \label{Xpde}
\begin{split}
\dt \phi_t =&\ \gD_{\til{g}(t), g_0} \phi_i(t) + \mathcal O \left(\del \phi
\right)\\
\phi_0 =&\ \id_M
\end{split}
\end{gather}
where $\gD_{\til{g}(t), g_0}$ is the harmonic map Laplacian taken with respect
to the metrics $\til{g}(t)$ and $g_0$.

We now proceed with the proof of uniqueness.  Let $\til{g}_1(t), \til{g}_2(t)$
be two
solutions to (\ref{flow}) with $\til{g}_1(0) = \til{g}_2(0) = g_0$.  Let
$\phi_i(t)$ be solutions to (\ref{Xpde}) with respect to $\til{g}_i$, which
exists in general because (\ref{Xpde}) is strictly parabolic and $M$ is compact.
 Now, pushing forward by these diffeomorphisms we observe that $g_i(t) :=
\left(\phi_i(t) \right)_* \til{g}_i(t)$ are both solutions of (\ref{gfflow}).
Since $g_1(0) = g_2(0)$ and solutions to (\ref{gfflow}) are unique, it follows
that $g_1(t) = g_2(t)$ as long as these metrics are defined.  But now one
observes that $\phi_1(t)$ and $\phi_2(t)$ are both solutions to the same ODE
(\ref{Xode}) with the same initial condition, and are therefore equal.  It
follows that $\til{g}_1(t) = \til{g}_2(t)$ as long as they are both defined and
the result follows.

Now consider the case where $J_0$ is integrable.  In this case one can consider
the flow of K\"ahler forms
\begin{align*}
\dt \omega =&\ -S + Q\\
\omega(0) =&\ \omega_0.
\end{align*}
Here $\omega_0$ is compatible with an integrable complex structure, and so this
equation is an example of \emph{Hermitian curvature flow}, studied in
\cite{ST1}.  In particular, there is a short time solution to this equation
which remains compatible with $J_0$.  If we let $J(t) = J_0$, it is easy to
verify that $(\omega(t), J(t))$ is a solution to (\ref{flow}), and more to the
point, the unique solution.  If the initial structure $(\omega, J)$ is K\"ahler
one can verify that the solution to (\ref{flow}) is the K\"ahler Ricci flow by a
similar argument.
\end{proof}

%%%%%%%%%%%%%%%%%%%%%%%%%%%%%%%%%%%
\section{Symplectic Curvature Flow} \label{symplecticste}
%%%%%%%%%%%%%%%%%%%%%%%%%%%%%%%%%%%

In this section we will motivate and investigate the equation (\ref{AKflow}).
The general philosophical starting point is clear: one would
like to define a flow of symplectic structures $\dt \omega = - P$, purely in
analogy with
K\"ahler Ricci flow.
However, $P$ is not a $(1,1)$-form, so $\omega$ would not stay compatible with
$J$, and then the definitions fall apart.  Thus one is naturally led to allowing
$J$
to flow as well.  Lemma \ref{formcompatibility} suggests \emph{part, but not
all} of what should appear in the evolution equation for $J$ (see Remark
\ref{comprmk}).  It is a small miracle that in the almost
K\"ahler setting there is a very natural choice for this component, i.e. the
endomorphism $\mathcal R$, which ends
up yielding a parabolic equation.  These considerations lead to the definition
of (\ref{AKflow}), which
we repeat here:
\begin{gather*}
\begin{split}
\dt \omega =&\ - { P}\\
\dt J =&\ - D^* D J - \mathcal N + \mathcal R\\
\omega(0) =&\ \omega_0\\
J(0) =&\ J_0.
\end{split}
\end{gather*}
Before proving Theorem \ref{AKflowthm} we give two equivalent formulations of
this flow, one putting it in the
framework of equation (\ref{almostcomplexflow}), the other realizing this system
as a flow coupling a parabolic for $J$ with Ricci flow.  The latter
formulation is the appropriate viewpoint to use to show the short time existence
of solutions to (\ref{AKflow}).

\begin{prop} \label{flowequiv} Let $(M^{2n}, \omega(t), J(t))$ be a
one-parameter family of almost K\"ahler structures solving (\ref{AKflow}).  Then
the family $(\omega(t), J(t))$ is a solution to
\begin{gather} \label{specialflow}
\begin{split}
\dt \omega =&\ - S + Q^1 + H\\
\dt J =&\ - \mathcal K + \mathcal H^1
\end{split}
\end{gather}
where $Q^1$ and $\mathcal H^1$ are defined in (\ref{Q1def}) and (\ref{QQ1def})
respectively, and $H$ is defined according to (\ref{Fdef}).  In particular,
(\ref{specialflow}) is a degenerate parabolic equation for almost
Hermitian pairs $(\omega, J)$ which preserves the almost K\"ahler condition.
\begin{proof} Let $(\omega, J)$ be an almost K\"ahler structure, and let
$\Omega$ denote the curvature of the Chern connection.  Let $\{e_i \}$ denote a
local orthonormal frame for $T^{1,0}(M)$.  First recall the Bianchi identity for
a connection $\N$:
\begin{align*}
\Sigma_{X,Y,Z} \left[ \Omega(X,Y)Z - T(T(X,Y),Z) - \N_X T(Y,Z) \right] = 0
\end{align*}
For our almost K\"ahler structure the torsion $T$ is completely determined by
the Nijenhuis tensor, which is a $(0,2)$ form with values in $(1,0)$ vectors.
Using this we compute an expression for the $(1,1)$ part of $P$.
\begin{gather} \label{PtoS}
\begin{split}
P(e_j, \bar{e}_k) =&\ \Omega(e_j, \bar{e}_k,e_i, \bar{e}_i)\\
=&\ \Omega(e_i, \bar{e}_k, e_j, \bar{e}_i) + \left< N(N(e_i, e_j), \bar{e}_k),
\bar{e}_i \right> + \left<\N_{\bar{e}_k} N(e_i, e_j), \bar{e}_i \right>\\
=&\ - \Omega(e_i, \bar{e}_k, \bar{e}_i, e_j) + \left< N(N(e_i, e_j), \bar{e}_k),
\bar{e}_i \right> + \left<\N_{\bar{e}_k} N(e_i, e_j), \bar{e}_i \right>\\
=&\ S(e_j, \bar{e}_k) - \left< N(N(\bar{e}_i, \bar{e}_k), e_i), e_j \right> -
\left< \N_{e_i} N(\bar{e}_k, \bar{e}_i), e_j \right>\\
&\ + \left< N(N(e_i, e_j), \bar{e}_k), \bar{e}_i \right> + \left<\N_{\bar{e}_k}
N(e_i, e_j), \bar{e}_i \right>.
\end{split}
\end{gather}
But since $N$ takes values in $(1,0)$ vectors and $\N$ is a Hermitian
connection, it follows that
\begin{align*}
\left< \N_{e_i} N(\bar{e}_k, \bar{e}_i), e_j \right> = \left< \N_{\bar{e}_k} N
(e_i, e_j), \bar{e}_i \right> = 0.
\end{align*}
It follows that
\begin{align*}
P^{1,1} =&\ S + Q^1
\end{align*}
where
\begin{align} \label{Q1def}
Q^1_{i \bj} =&\ \omega^{k \bl} \left( g_{m \bl} N_{k i}^{\bp} N_{\bp \bj}^{m} -
g_{m \bj} N_{\bl \bj}^{p} N_{p k}^{\bm} \right).
\end{align}
Next we examine the evolution equation for $J$.  Choose normal coordinates for
the associated metric at a point $p$.  Then, including the precise lower order
terms in (\ref{Kcalc}), we see that
\begin{align*}
\omega^{kl} \N_k N_{l j}^i =&\ - g^{kl} \del_k \del_l J_j^i + \omega^{kl} \left(
J_q^i \del_k \del_j J_l^q - J_j^q \del_k \del_q J_l^i \right)\\
&\ + \omega^{kl} \left( D_k J_l^p D_p J_j^i - D_k J_j^p D_p J_l^i - D_k J_p^i
D_l J_j^p + D_k J_p^i D_j J_l^i \right)\\
&\ + \frac{1}{2} \omega^{kl} \left( N_{k p}^i N_{l j}^p - N_{k l}^p N_{p j}^i -
N_{k j}^p N_{l p}^i \right).
\end{align*}
Furthermore, by a direct calculation in
normal coordinates at $p$ one has
\begin{align*}
( - D^* D J + \mathcal N + \mathcal R)_{j}^i =&\ g^{kl} \del_k \del_l J_j^i +
g^{kl} J_j^p \del_p \gG_{kl}^i - g^{kl} J_p^i \del_j \gG_{kl}^p + \mathcal
N_j^i.
\end{align*}
Furthermore, calculating as in (\ref{Xcalc}), (\ref{Xcalc2}), again using the
normal coordinates,
\begin{align*}
g^{kl} J_j^p \del_p \gG_{kl}^i =&\ J_j^p \del_p \left( g^{kl} \gG_{kl}^i
\right)\\
=&\ J_j^p \del_p \left( g^{kl} (\gG_C)_{kl}^i + \omega^{kl} (\gG_C)_{kl}^q J_q^i
\right)\\
=&\ J_j^p \del_p \left( \omega^{kl} \del_k J_l^i - \gw^{kl} \N_k J_l^i \right)\\
=&\ J_j^p \omega^{kl} \del_p \del_k J_l^i - J_j^p \omega^{kr} J_q^l D_p
J_{r}^q D_k J_l^i.
\end{align*}
Likewise
\begin{align*}
g^{kl} J_p^i \del_j \gG_{kl}^p = \omega^{kl} J_p^i \del_j \del_k J_l^p - J_p^i
\omega^{kr} J_q^l D_j J_r^q D_k J_l^p.
\end{align*}
Combining these calculations yields
\begin{align*}
\left( - D^* D J + \mathcal N + \mathcal R \right) =&\ - \mathcal K + \mathcal
H^1
\end{align*}
where
\begin{gather} \label{QQ1def}
\begin{split}
\left(\mathcal H^1 \right)_j^i =&\ \omega^{kl} \left( D_k J_l^p D_p J_j^i - D_k
J_j^p D_p J_l^i - D_k J_p^i
D_l J_j^p + D_k J_p^i D_j J_l^i \right)\\
&\ + \frac{1}{2} \omega^{kl} \left( N_{k p}^i N_{l j}^p - N_{k l}^p N_{p j}^i -
N_{k j}^p N_{l p}^i \right)\\
&\ - J_j^p \omega^{kr} J_q^l D_p J_{r}^q D_k J_l^i + J_p^i
\omega^{kr} J_q^l D_j J_r^q D_k J_l^p + g^{ik} g^{pq} \omega_{rs} D_p J_j^r D_q
J_k^s.
\end{split}
\end{gather}
Finally, it is clear by construction that if we define $H$ so that (\ref{Fdef})
holds, it must equal $- P^{2,0 + 0,2}$.  It follows that a solution to
(\ref{AKflow}) is a solution to (\ref{specialflow}), and the proposition
follows.  The final statement of the proposition will follow once existence and
uniqueness of solutions to (\ref{AKflow}) is established.
\end{proof}
\end{prop}

\begin{prop} \label{Bfieldprop1} Let $(M^{2n}, \omega(t), J(t))$ be a
one-parameter family of almost
K\"ahler structures solving (\ref{AKflow}).  Then the associated Riemannian
metric $g(t)$ satisfies
\begin{align} \label{sympmetricev}
\dt g =&\ -2 \Ric + \frac{1}{2} B^1 - B^2
\end{align}
where
\begin{align*}
B^i(\cdot, \cdot) = N^i(\cdot, J \cdot).
\end{align*}
In coordinates one has
\begin{align*}
B^1_{ij} =&\ g^{kl} g_{mn} D_i J_k^m D_j J_l^n,\\
B^2_{ij} =&\ g^{kl} g_{mn} D_k J_i^m D_l J_j^n.
\end{align*}

\begin{proof} We begin with a general calculation using the notation of Lemma
\ref{formcompatibility}.  Specifically, we have
\begin{align*}
\dt g(\cdot, \cdot) =&\ \dt \left[\omega (\cdot, J \cdot ) \right]\\
=&\ \left[ \phi (\cdot, J \cdot ) + \psi(\cdot, J \cdot) + \omega(\cdot, K
\cdot) \right].
\end{align*}
Let us compute these three terms separately.  First of all, since $d \omega = 0$
it follows from
(\ref{P20part}) and (\ref{Pformula}) that
{
\begin{align*}
\phi(\cdot, J \cdot) =&\ - P^{1,1} (\cdot, J \cdot)\\
=&\  P^{2,0 + 0,2} (\cdot, J
\cdot) - P(\cdot, J \cdot)\\
=&\ \left( - 2 \rho + \frac{1}{2} N^1 {- N^2}
\right) \left(\cdot, J \cdot \right)\\
=&\ -2 \Ric^J (\cdot, \cdot) + \frac{1}{2} B^1 (\cdot, \cdot) - B^2(\cdot,
\cdot).
\end{align*}}
Now observe that
\begin{align*}
\psi \left( \cdot, J \cdot \right) =&\ {- P^{2,0 + 0,2} \left(\cdot, J \cdot
\right)}.
\end{align*}
Next consider
\begin{align*}
\omega(\cdot, K \cdot)_{ij} = \omega_{ik} K^k_j = \omega_{ik} \left( g^{kl}
\left( { - P^{2,0 + 0,2}_{jl}} \right) + J_j^l \Rc_l^k - \Rc_j^l J_l^k  \right).
\end{align*}
The first term simplifies to
{
\begin{align*}
- \omega_{ik} g^{kl} P_{jl}^{2,0 + 0,2} =&\ - J_i^p g_{pk} g^{kl} P_{jl}^{2,0
+ 0,2}\\
=&\ -J_i^l P^{2,0 + 0,2}_{jl}\\
=&\ J_i^l \left( J_j^m J_l^p P^{2,0 + 0,2}_{mp} \right)\\
=&\ - J_j^m P_{m i}^{2,0 + 0,2}\\
=&\ P^{2,0 + 0,2} \left(\cdot, J \cdot \right)_{ij}.
\end{align*}}
Next we calculate
\begin{align*}
\omega_{ik} \left( \Rc_j^l J_l^k - J_j^l \Rc_l^k \right) =&\ J_i^p g_{pk} \left(
J_j^l \Rc_l^k - \Rc_j^l J_l^k \right)\\
=&\ \left(J^* \Ric - \Rc \right)_{ij}\\
=&\ - 2 \left( \Rc - \Ric^J \right)_{ij}.
\end{align*}
Combining the above calculations, the result follows.
\end{proof}
\end{prop}

\noindent Now we give the proof of Theorem \ref{AKflowthm}.

\begin{proof} As in the proof of Theorem \ref{almosthermitianste}, we must work
in a more general setting to construct solutions to (\ref{AKflow}).  It will be
most convenient in this case to work with a coupled system of Riemannian metric
and almost complex structure.  The main goal will be to construct a solution to
\begin{gather} \label{riemsympflow}
\begin{split}
\dt g =&\ -2 \Ric + \left(\frac{1}{2} N^1(\cdot, J \cdot) - N^2 (\cdot, J \cdot)
\right)^J\\
\dt J =&\ - D^* D J - \mathcal N + \mathcal R\\
g(0) =&\ g_0\\
J(0) =&\ J_0
\end{split}
\end{gather}
for a compatible initial condition $(g_0, J_0)$.  We will accomplish this in the
same two step process employed in the proof of Theorem \ref{almosthermitianste},
first defining a gauge-fixed flow, then defining a flow for arbitrary pairs $(g,
J)$ which preserves compatibility of the initial condition. In particular, let
$X$ be defined as in (\ref{Xdef}), and consider the gauge fixed flow
\begin{gather} \label{gfriemsympflow}
\begin{split}
\dt g =&\ -2 \Ric + \left( \frac{1}{2} N^1(\cdot, J \cdot) - N^2 (\cdot,
J \cdot) \right)^{J} + L_X g =: \mathcal D_1(g, J)\\
\dt J =&\ - D^* D J - \mathcal N + \mathcal R + L_X J =: \mathcal D_2(g, J).
\end{split}
\end{gather}
Also as in Theorem \ref{almosthermitianste} we must modify the flow and pull
back to a linear space to construct the solution to (\ref{gfriemsympflow}).
Using the notation of the proof of Theorem \ref{almosthermitianste}, let
\begin{gather} \label{symplpfloc10}
\begin{split}
\dt g =&\ \mathcal D_1(g^{\pi E}, \pi E) - D^*_{g_0} D_{g_0} g^{-\pi E} =:
\widetilde{\mathcal D_1}(g, E)\\
\dt E =&\ (D \pi_{\pi E}^{-1}) \left( \mathcal D_2(g^{\pi E}, \pi E) \right) =:
\widetilde{\mathcal D_2}(g, E)\\
g(0) =&\ g_0\\
E(0) =&\ 0.
\end{split}
\end{gather}
First note that in fact the evolution for $E$ is well defined.  Indeed, an
endomorphism
of the form $g^{-1} \psi$ where $\psi \in \Lambda^{(2,0) + (0,2)}$ automatically
satisfies $K J + J K = 0$.  Therefore by Lemma \ref{P2formula} the tensor
$-D^* D J + \mathcal N$ is already $J$-skew, as of course is $\mathcal R$.
Therefore $\mathcal D_2(g^{\pi E}, \pi E) \in \mathbb T \mathbb J_{\pi E}$, and
so $\widetilde{\mathcal D}_2$ is well defined.

We compute the linearization of (\ref{symplpfloc10}) at $t = 0$.  Here the small
miracle of equation (\ref{Xform}) is highly relevant, i.e. that
the vector field $X$ generates the same one parameter family of diffeomorphisms
which appears in proving short time existence of Ricci flow.  It follows from a
well-known calculation (see \cite{Chow} pg. 114) that
\begin{align*}
\mathcal L_g (\left(-2 \Rc + L_X g \right) (h) =&\ \gD_{L} h + \mbox{ lower
order terms.}
\end{align*}
Since the terms $N^1$ and $N^2$ are both first order in both $g$ and $J$, we
conclude that
\begin{align*}
\sigma \left[ \widehat{\mathcal L_{g} \widetilde{\mathcal D}_1} \right](h)_{ij}
=&\
\brs{\xi}^2 h^{\pi E}_{ij} + \brs{\xi}^2 h^{-\pi E}_{ij} = \brs{\xi}^2 h_{ij}\\
\sigma \left[ \widehat{\mathcal L_{J} \widetilde{\mathcal D}_1} \right](K)_{ij}
=&\ 0.
\end{align*}
Next we compute the linearization of $\widetilde{\mathcal D}_2$ at $t = 0$.  We
do this by first computing a coordinate formula for $\mathcal D_2$ applied to a
compatible pair $(g, J)$ in stages.  First
of all
we have
\begin{gather}
\begin{split}
\left[- D^* D J\right]_k^l =&\ g^{i j} \del_i \left[ D J
\right]_{j k}^l + \mathcal O(\del J, \del g)\\
=&\ g^{i j} \del_i \left[ \del_j J_k^l - \gG_{j k}^p J_p^l + \gG_{j
p}^l J_k^p \right] + \mathcal O(\del J, \del g)\\
=&\ g^{i j} \left[\del_i \del_j J_k^l - \del_i \gG_{j k}^p J_p^l +
\del_i \gG_{j p}^l J_k^p \right] + \mathcal O(\del J, \del g)
\end{split}
\end{gather}
where here $\gG$ denotes the Levi Civita connection.  Next, using Lemma
\ref{coordianteLieJ} and (\ref{Xform}) we compute
\begin{align*}
\left[\mathcal L_X J \right]_k^l =&\ J_p^l \del_k \left(
g^{i j} \gG_{i j}^p \right) - J_k^p \del_p \left( g^{i j} \gG_{i
j}^l \right) + \mathcal O(\del J, \del g)\\
=&\ J_p^l g^{ij} \del_k \gG_{i j}^p - J_k^p g^{i j}
\del_p \gG_{i j}^l + \mathcal O(\del J, \del g).
\end{align*}
Combining the above calculations we observe
\begin{gather} \label{D2calc}
\begin{split}
\mathcal D_2(g, J)_k^l =&\ g^{i j} \left[\del_i \del_j J_k^l - \del_i
\gG_{j k}^p J_p^l +
\del_i \gG_{j p}^l J_k^p \right] + J_p^l g^{ij} \del_k \gG_{i j}^p - J_k^p g^{i
j}
\del_p \gG_{i j}^l\\
&\ + J_k^p \Rc_p^l - \Rc_k^p J_p^l + \mathcal O(\del J, \del g)\\
=&\ g^{ij} \del_i \del_j J_k^l + g^{ij} J_k^p \left[\del_i \gG_{j p}^l - \del_p
\gG_{i j}^l + \Rm_{pij}^l \right]\\
&\ + g^{ij} J_p^l \left[ \del_k \gG_{i j}^p - \del_i \gG_{j k}^p - \Rm_{kij}^p
\right] + \mathcal O(\del J, \del g)\\
=&\ g^{ij} \del_i \del_j J_k^l + \mathcal O(\del J, \del g).
\end{split}
\end{gather}
It follows that
\begin{align*}
\sigma \left[ \widehat{\mathcal L_{g} \mathcal D_2} \right](h)_{i}^j =&\
0\\
\sigma \left[ \widehat{\mathcal L_{J} \mathcal D_2} \right](K)_{i}^j =&\
\brs{\xi}^2 K_i^j.
\end{align*}
Using that $D \pi_0 = \Id_{\mathbb T \mathbb J_{J_0}}$ it follows that the same
formulas hold for the linearization of $\widetilde{\mathcal D}_2$.  We conclude
\begin{align*}
\gs \left[ \widehat{\mathcal L \widetilde{\mathcal D}} \right](h, K) =&\
\left(
\begin{matrix}
I & 0\\
0 & I
\end{matrix}
\right) \left(
\begin{matrix}
h\\
K
\end{matrix}
\right).
\end{align*}
It follows from standard parabolic theory that,
starting from a compatible pair $(g_0, J_0)$, there is a unique short time
solution to (\ref{symplpfloc10}).

Now let $J = \pi E$.  We want to show that the pair $(g, J)$ is a solution to
(\ref{gfriemsympflow}).  By a straightforward computation using that $g(D^* D J
- \mathcal N, \cdot) \in \Lambda^{(2,0) + (0,2)}$, one
can show that the differential operators $\mathcal D_1$ and $\mathcal D_2$
satisfy, for a compatible pair $(g, J)$,
\begin{align*}
\mathcal D_1(g^J, J)^{-J} + g(J \cdot, \mathcal D_2(g^J, J) \cdot) +
g(\mathcal D_2(g^J, J) \cdot, J \cdot) = 0.
\end{align*}
It follows that a solution to (\ref{symplpfloc10}) satisfies
\begin{align*}
\dt g^{-J} =&\ D^*_{g_0} D_{g_0} g^{-J},
\end{align*}
and at this point one can follow the estimate of (\ref{genpfloc10}) and apply
the maximum principle to obtain that compatibility of the pair $(g_0, J_0)$ is
preserved along a solution to (\ref{symplpfloc10}).  It follows that $(g(t),
J(t))$ is a solution of (\ref{gfriemsympflow}).
If $\phi_t$ denotes the one parameter family of diffeomorphisms generated by
$-X$ as in (\ref{Xode}), then $(\phi_t^* g(t), \phi_t^* J(t))$ is a solution to
(\ref{riemsympflow}).  Furthermore, the proof of uniqueness of solutions to
(\ref{riemsympflow}) follows exactly the same lines as the uniqueness part of
Theorem \ref{almosthermitianste}, since it is the same vector field $X$ we are
using in the gauge-fixing technique.

Next we must show that the symplectic condition $d\omega = 0$ is preserved along
our solution to (\ref{riemsympflow}).  We want to reverse the steps of
Proposition \ref{Bfieldprop1}, except this time we need to use more general
formula since our K\"ahler form $\omega$ is not a priori symplectic.  We start
with
\begin{align*}
\dt \omega =&\ \dt g (J \cdot, \cdot)\\
=&\ - 2 \Rc (J \cdot, \cdot) + \frac{1}{2} N^1(\cdot, \cdot)^{J} - N^2(\cdot,
\cdot)^J + g(
- D^* D J {+ \mathcal N} + \mathcal R, \cdot)
\end{align*}
As $\mathcal R$ is the $J$-skew part of the Ricci tensor, one has
\begin{align*}
-2 \Rc(J \cdot, \cdot) + g(\mathcal R \cdot, \cdot) = -2 \rho (\cdot, \cdot).
\end{align*}
Furthermore, as the calculation of Lemma \ref{P2formula} shows, the tensor $N^2$
is in general already $J$-symmetric, and $\mathcal N$ is just $g^{-1} N^2$, so
it
follows that
\begin{align*}
- N^2(\cdot, \cdot)^J + g(\mathcal N, \cdot)= 0.
\end{align*}
Combining these calculations yields
{
\begin{align*}
\dt \omega =&\ - 2 \rho + \left( \frac{1}{2} N^1 \right)^J - D^* D \omega\\
=&\ - \rho^* + \left( \frac{1}{2} N^1 \right)^J - \gD_d \omega\\
=&\ - P + \left( \frac{1}{2} N^1 \right)^{-J} + \frac{1}{2} W - \gD_d \omega.
\end{align*}}
Now note that it follows from Lemma \ref{generalJlemma} that $W = d \omega *
DJ$, hence
\begin{align*}
\dt d \omega =&\ - d \gD_d \omega + A * d \omega + B * \N d \omega\\
=&\ - \gD_d d \omega + A * d \omega + B * D d \omega.
\end{align*}
for some tensor quantities $A$ and $B$.  It follows using the B\"ochner formula
that
\begin{align*}
\dt \brs{d \omega}^2_g =&\ - 2 \left< \gD_d d \omega, d \omega \right> + 2
\left< A * d \omega + B * D d \omega, d \omega \right>\\
=&\ - 2 \left< D^* D \omega + \Rm * d \omega, d \omega \right> + 2 \left< A * d
\omega + B * D d \omega, d \omega \right>\\
\leq&\ \gD \brs{ d \omega}^2 - 2 \brs{D d \omega}^2 + C \brs{d
\omega}^2 + 2 B \brs{D d \omega} \brs{d \omega}\\
\leq&\ \gD \brs{ d \omega}^2 - 2 \brs{D d \omega}^2 + C \brs{d
\omega}^2 + \left(C \ge \brs{D d \omega}^2 + \frac{C}{\ge} \brs{d \omega}^2
\right)\\
\leq&\ \gD \brs{d \omega}^2 - \brs{D d \omega}^2 + C \brs{d
\omega}^2
\end{align*}
where the last line follows by choosing $\ge$ small with respect to bounds on
$A$ and $B$.  Applying the maximum principle to $e^{-C t} \brs{d \omega}_g^2$,
it follows that $d \omega \equiv 0$ is
preserved along the solution to (\ref{riemsympflow}).

The final step is to show that if the initial structure is K\"ahler then the
solution reduces to K\"ahler Ricci flow.  Given $(\omega_0, J_0)$ K\"ahler, we
can construct the solution to (\ref{AKflow}) by the above proof.  However we can
also construct the solution to K\"ahler Ricci flow with the initial condition
$\omega_0$.  Since the K\"ahler condition is preserved, one observes that in
fact $(\omega(t), J_0)$ is a solution to (\ref{AKflow}).  Since solutions to
(\ref{AKflow}) are unique, it follows that our given solution to (\ref{AKflow})
must be the same as the solution to K\"ahler Ricci flow.
\end{proof}

\begin{rmk} We close this section with an important remark regarding the
uniqueness of equations satisfying the results of Theorem \ref{AKflowthm}.  In
principle, there is a natural \emph{family} of parabolic flows which
preserve the almost K\"ahler condition for the pair $(\omega, J)$.  To see this
observe that if $\mathcal H$ denotes any endomorphism of the tangent bundle
consisting only of first order terms in $\omega$ and $J$, which further
satisfies
\begin{gather} \label{kernel}
\begin{split}
\omega(\mathcal H, J) + \omega(J, \mathcal H) = 0, \qquad J \mathcal H +
\mathcal H J = 0
\end{split}
\end{gather}
then the proof of Theorem \ref{AKflowthm} carries through for the evolution
equation
\begin{gather} \label{generalAKflow}
\begin{split}
\dt \omega =&\ - P\\
\dt J =&\ - D^* D J + \mathcal N + \mathcal R + \mathcal H\\
\omega(0) =&\ \omega_0\\
J(0) =&\ J_0
\end{split}
\end{gather}
\noindent Indeed, it follows from (\ref{kernel}), the prior calculations and
Lemma
\ref{formcompatibility} that compatibility of the pair will still be preserved
under this flow.  One can think of this also as adding a certain $(2,0) + (0,2)$
tensor
to the evolution of both $g$ and $J$, which cancels out and thus does not appear
in
the evolution of $\omega$.

Furthermore, since the term $\mathcal H$ is first
order in
$\omega$ and $J$ the discussion of parabolicity and short time existence is not
affected.  Observe that since we are staying within the class of almost K\"ahler
structures, the only first order invariant of the pair $(\omega, J)$ is the
Nijenhuis tensor.  Therefore in principle we could let $\mathcal H$ denote any
expression in the Nijenhuis tensor which satisfies (\ref{kernel}).  To produce
an evolution equation with a natural scaling property, it is most relevant to
consider endomorphisms $\mathcal H$ which are \emph{quadratic}.  We codify this
discussion with the following proposition.

\begin{prop} Let $(M^{2n}, \omega_0, J_0)$ denote an almost K\"ahler structure.
Suppose $\mathcal H \in \End(TM)$ is a quadratic expression in the Nijenhuis
tensor satisfying (\ref{kernel}).  Then there exists a unique short time
solution to (\ref{generalAKflow}) with initial condition $(\omega_0, J_0)$.
\end{prop}

\end{rmk}

%%%%%%%%%%%%%%%%%%%%%%%%%%%%%%%%%%%%%%%
\section{Curvature Evolution Equations} \label{curvatureev}
%%%%%%%%%%%%%%%%%%%%%%%%%%%%%%%%%%%%%%%

In this section we derive evolution equations for the curvature of the Chern
connection, the Nijenhuis tensor, and their derivatives for solutions to
(\ref{almostcomplexflow}).  By Proposition \ref{flowequiv}, these general
equations hold for solutions to the symplectic curvature flow as well.  As one
would expect from the calculation of the symbol of (\ref{almostcomplexflow}) in
Theorem \ref{almosthermitianste}, the system of equations is upper triangular in
the appropriate sense.

\subsection{Evolution equations for almost Hermitian curvature flow}

First we derive the evolution of the Nijenhuis tensor.

\begin{prop} \label{Nevolution} Let $(M^{2n}, \omega(t), J(t))$ be a solution to
(\ref{almostcomplexflow}).  Then
\begin{align} \label{Nijenhuisev}
\dt N =&\ \gD N + \Omega * T + \N T * T + T^{*3}.
\end{align}
\begin{proof} Choose normal coordinates for the induced metric $g =
\omega(J\cdot, \cdot)$ at a fixed point in space and time.  Note that this has
the effect that any first coordinate derivative of $\omega$ or $J$ can be
expressed in terms of the torsion of the Chern connection $T$.  Furthermore, any
second coordinates derivative of $J$ can be expressed using the curvature,
torsion, and first derivative of torsion.  Starting from (\ref{Nijencoords}) we
compute
\begin{align*}
\dt N_{jk}^i =&\ \left( \dot{J} * \del J \right)_{jk}^i + J_j^p \del_p
\dot{J}_{k}^i - J_k^p \del_p \dot{J}_j^i - J_p^i \del_j \dot{J}_k^p + J_p^i
\del_k \dot{J}_j^p.
\end{align*}
Since $\dot{J} = - \mathcal K + \mathcal H$ where $\mathcal H$ is quadratic in
the torsion of $\N$, it immediately follows that
\begin{align*}
\dot{J} * \del J = T^{*3}
\end{align*}
where $T$ denotes the full torsion tensor of the Chern connection.  Of course in
the symplectic setting this only depends on $N$.  Furthermore, it is clear that
\begin{align*}
\del \mathcal H = \N T * T + T^{*3}.
\end{align*}
Therefore it remains to calculate the highest order term
\begin{align*}
\mathcal W_{jk}^i :=&\ - J_j^p \del_p \mathcal K_{k}^i + J_k^p \del_p \mathcal
K_j^i + J_p^i \del_j \mathcal K_k^p - J_p^i \del_k \mathcal K_j^p\\
=&\ - J_j^p \del_p \left( \omega^{rs} \N_r N_{s k}^i \right) + J_k^p \del_p
\left( \omega^{rs} \N_r N_{s j}^i \right) + J_p^i \del_j \left( \omega^{rs} \N_r
N_{s k}^p \right) - J_p^i \del_k \left( \omega^{rs} \N_r N_{s j}^p \right)\\
=&\ - J_j^p \del_p \left( \omega^{rs} \del_r N_{s k}^i + T^{*2} \right)+ J_k^p
\del_p \left( \omega^{rs} \del_r N_{s j}^i + T^{*2} \right)\\
&\ + J_p^i \del_j \left( \omega^{rs} \del_r N_{s k}^p + T^{*2} \right) - J_p^i
\del_k \left( \omega^{rs} \del_r N_{s j}^p + T^{*2} \right)\\
=&\ - J_j^p \omega^{rs} \del_p \del_r N_{sk}^i + J_k^p \omega^{rs} \del_p \del_r
N_{s j}^i + J_p^i \omega^{rs} \del_j \del_r N_{sk}^p - J_p^i \omega^{rs} \del_k
\del_r N_{sj}^p + \N T * T.
\end{align*}
Now plugging in (\ref{Nijencoords}) again we conclude
\begin{align*}
\mathcal W_{jk}^i =&\ - J_j^p \omega^{rs} \del_p \del_r \left( J_s^t \del_t
J_k^i - J_k^t \del_t J_s^i - J_t^i \del_s J_k^t + J_t^i \del_k J_s^t \right)\\
&\ + J_k^p \omega^{rs} \del_p \del_r \left( J_s^t \del_t J_j^i - J_j^t \del_t
J_s^i - J_t^i \del_s J_j^t + J_t^i \del_j J_s^t\right)\\
&\ + J_p^i \omega^{rs} \del_j \del_r \left(J_s^t \del_t J_k^p - J_k^t \del_t
J_s^p - J_t^p \del_s J_k^t + J_t^p \del_k J_s^t \right)\\
&\ - J_p^i \omega^{rs} \del_k \del_r \left(J_s^t \del_t J_j^p - J_j^t \del_t
J_s^p - J_t^p \del_s J_j^t + J_t^p \del_j J_s^t\right) + \N T * T\\
=&\ - J_j^p \omega^{rs} \left( J_s^t \del_p \del_r \del_t J_k^i - J_k^t \del_p
\del_r \del_t J_s^i - J_t^i \del_p \del_r \del_s J_k^t + J_t^i \del_p \del_r
\del_k J_s^t \right)\\
&\ + J_k^p \omega^{rs} \left( J_s^t \del_p \del_r \del_t J_j^i - J_j^t \del_p
\del_r \del_t J_s^i - J_t^i \del_p \del_r \del_s J_j^t + J_t^i \del_p \del_r
\del_j J_s^t \right)\\
&\ + J_p^i \omega^{rs} \left( J_s^t \del_j \del_r \del_t J_k^p - J_k^t \del_j
\del_r \del_t J_s^p - J_t^p \del_j \del_r \del_s J_k^t + J_t^p \del_j \del_r
\del_k J_s^t \right)\\
&\ - J_p^i \omega^{rs} \left( J_s^t \del_k \del_r \del_t J_j^p - J_j^t \del_k
\del_r \del_t J_s^p - J_t^p \del_k \del_r \del_s J_j^t + J_t^p \del_k \del_r
\del_j J_s^t \right)\\
&\ + \Omega * T + \N T * T.
\end{align*}
Let us label the sixteen third derivative terms above as $I - XVI$ in Roman
numerals in the order in which they appear.  Some cancellations are apparent,
namely, using that $\omega$ is skew symmetric and coordinate derivatives are
symmetric, it
follows that $III = 0, VII = 0, XI = 0, XV = 0$.  Also, one observes that $XII +
XVI = 0$, $II + VI = 0$, $IV + XIV = 0$ and $VIII + X = 0$.  Therefore
\begin{align*}
\mathcal W_{jk}^i =&\ g^{rt} \del_r \del_t \left( J_j^p \del_p J_k^i - J_k^p
\del_p J_j^i - J_p^i \del_j J_k^p + J_p^i \del_k J_j^p \right) + \Omega * T + \N
T * T\\
=&\ g^{r t} \del_r \del_t N_{jk}^i + \Omega * T + \N T * T\\
=&\ \gD N_{j k}^i + \Omega * T + \N T * T.
\end{align*}
The result follows.
\end{proof}
\end{prop}

\begin{prop} \label{domevolution} Let $(M^{2n}, \omega(t), J(t))$ be a solution
to
(\ref{almostcomplexflow}).  Then
\begin{align} \label{domev}
\dt \left( d \omega \right)^+ =&\ \gD \left( d \omega \right)^+ + \N^2 N +
\Omega * T + \N T * T + T^{*3}.
\end{align}
\begin{proof} We will use a specialized frame to make this calculation.  We
choose normal coordinates for the induced metric at a given point as in
Proposition \ref{Nevolution}.  Let $\{\frac{\del}{\del x_i} \}_{i = 1, \dots,
2n}$ be
the associated coordinate vector fields .  Now let
\begin{gather*}
e_j := \frac{\del}{\del x_j} - i J(0) \left(\frac{\del}{\del x_j} \right).
\end{gather*}
Clearly $\{e_j\}_{j = 1}^{2n}$ contains a spanning set for $T^{1,0}(M)$, and by
relabeling we may assume that $\{e_i \}_{j = 1}^n$ is a basis for $T^{1,0}(M)$.
 One has the
corresponding basis for $T^{0,1}(M)$ given by
\begin{gather*}
\bar{e}_j = \frac{\del}{\del x_j} + i J(0) \left(\frac{\del}{\del x_j}
\right).
\end{gather*}
Now note that
\begin{align*}
\left( d \omega \right)^+ =&\ \left( d \omega + d \omega (\cdot, J \cdot, J
\cdot) + d \omega(J \cdot, \cdot, J \cdot) + d \omega(J \cdot, J \cdot, Z)
\right).
\end{align*}
It follows that
\begin{align*}
\dt \left( d \omega \right)^+ = \dot{J} * d \omega + \left( \dt d \omega
\right)^{(2,1) + (1,2)}=&\ \left( \dt d \omega \right)^{(2,1) + (1,2)} + \N T *
T.
\end{align*}
We now compute
\begin{align*}
\left( \dt d \omega \right)_{i j \bk} =&\ \del_i \dot{\omega}_{j \bk} - \del_j
\dot{\omega}_{i \bk} - \del_{\bk} \dot{\omega}_{i j}\\
=&\ \N_i \left( - S + Q + H \right)_{j \bk} - \N_j \left( - S + Q + H
\right)_{i \bk}\\
&\ - \N_{\bk} \left( - S + Q + H \right)_{ij} + T * \Omega + T^{*3} + T * \N
T\\
=&\ - \N_i S_{j \bk} + 2 \N_j S_{i \bk} + \N^2 N + T * \Omega + T^{*3} + T *
\N T
\end{align*}
where the last line follows using that $S \in \Lambda^{1,1}$, $Q$ is quadratic
in the torsion, and $H$ depends on a quadratic term in torsion and the
derivative of the Nijenhuis tensor.  Now we apply the Bianchi identity to
conclude
\begin{align*}
\N_j S_{i \bk} - \N_i S_{j \bk} =&\ \omega^{m\bn} \left[ \N_j \Omega_{m \bn i
\bk} - \N_i \Omega_{m \bn j \bk} \right]\\
=&\ \omega^{m \bn} \left[ \N_j \left( \Omega_{i \bn m \bk} + \N_{\bn} T_{i m
\bk} \right) - \N_i \left( \Omega_{j \bn m \bk} + \N_{\bn} T_{j m \bk} \right)
\right].
\end{align*}
Next we apply the differential Bianchi identity and simplify
\begin{align*}
\omega^{m \bn} \left(\N_j \Omega_{i \bn m \bk} - \N_i \Omega_{j \bn m \bk}
\right) =&\ \N_{\bn} \Omega_{i j m \bk} + T * \Omega\\
=&\ \N^2 N + \N T * T + T * \Omega
\end{align*}
where the last line follows because the $(2,0) + (0,2)$ component of the Chern
curvature only depends on a quadratic expression in the torsion and one
derivative of the Nijenhuis tensor.  We commute derivatives and apply the
Bianchi identity a final time to conclude
\begin{align*}
\omega^{m \bn} \left( \N_j \N_{\bn} T_{i m \bk} - \N_i \N_{\bn} T_{j m \bk}
\right) =&\ \omega^{m \bn} \left( \N_{\bn} \left( \N_j T_{i m \bk} - \N_i T_{j m
\bk} \right) \right) + \Omega * T\\
=&\ \omega^{m \bn} \left( \N_{\bn} \left( \N_m T_{i j \bk} + T^{*2} \right)
\right) + \Omega * T\\
=&\ \omega^{m \bn} \N_{\bn} \N_m T_{i j \bk} + \N T * T + \Omega * T\\
=&\ \gD T_{i j \bk} + \N T * T + \Omega * T.
\end{align*}
But, since $(d \omega)^-$ can be expressed using the Nijenhuis tensor, we
conclude from (\ref{generalconnection}) that $\gD T = \gD \left(d \omega
\right)^+ + \N^2 N$.  The result follows.
\end{proof}
\end{prop}

\begin{prop} \label{Tevolution} Let $(M^{2n}, \omega(t), J(t))$ be a solution to
(\ref{almostcomplexflow}).  Then
\begin{align} \label{Torsionev}
\dt T =&\ \gD T+ \N^2 N + \Omega * T + \N T * T + T^{*3}.
\end{align}
\begin{proof} Since by (\ref{generalconnection}) the torsion is determined
algebraically by the Nijenhuis tensor and $(d \omega )^+$, the result follows
immediately from Propositions \ref{Nijenhuisev} and \ref{domevolution}.
\end{proof}
\end{prop}

Next we will compute the evolution of the Riemannian curvature tensor.  To do
this we will first observe a general formula for $S$ in terms of the Ricci
tensor.

\begin{lemma} \label{generalS} Let $(M^{2n}, \omega, J)$ be an almost Hermitian
manifold.  Then
\begin{align*}
S(J \cdot, \cdot) =&\ 2 \Rc + D T + T^{*2}.
\end{align*}
\begin{proof} This is an immediate consequence of various lemmas above.
\end{proof}
\end{lemma}

\begin{prop} Let $(M^{2n}, \omega(t), J(t))$ be a solution to
(\ref{almostcomplexflow}).  Then
\begin{align} \label{curvev}
\dt \Rm =&\ \gD \Rm + \Rm^{*2} + \Rm * T^{*2} + \Rm * D T + T * D^2 T + D T * D
T + D^3 T.
\end{align}
\begin{proof} Recall the general variational formula for the Riemannian
curvature tensor.  If $\dt{g} = h$ then
\begin{gather} \label{curvvar}
\begin{split}
\dt \Rm_{ijkl} =&\ \frac{1}{2} \left( D_i D_k h_{j l} - D_i D_l h_{j k} - D_j
D_k h_{i l} + D_j D_l h_{i k} \right)\\
&\ + \frac{1}{2} \left(R_{i j k}^p h_{p l} - R_{i j l}^p h_{p k} \right).
\end{split}
\end{gather}
First recall that for variation by $-2 \Rc$, one has
\begin{align*}
\dt \Rm =&\ \gD \Rm + \Rm^{*2}.
\end{align*}
The precise form of the equation appears in the next subsection on the
symplectic flow.  Using Lemma \ref{generalS}, the result follows.
\end{proof}
\end{prop}

\begin{thm} \label{ACFderevolution} Let $(M^{2n}, \omega(t), J(t))$ be a
solution to
(\ref{almostcomplexflow}).  Then
\begin{align*}
\dt D^k N =&\ \gD D^k N + \sum_{i = 0}^{k-1} D^i \Rm * D^{k-i} T + \sum_{i =
1}^{k+1} D^i T * D^{k+1-i} T\\
&\ + \sum_{i = 0}^{k} \sum_{j = 0}^i D^j T * D^{i-j}T * D^{k-i} T,\\
\dt D^k T =&\ \gD D^k T + \sum_{i = 0}^{k-1} D^i \Rm * D^{k-i} T + \sum_{i =
1}^{k+1} D^i T * D^{k+1-i} T\\
&\ + \sum_{i = 0}^{k} \sum_{j = 0}^i D^j T * D^{i-j}T * D^{k-i} T + D^{k+2} N,\\
\dt D^k \Rm =&\ \gD D^k \Rm + \sum_{i = 0}^k D^i \Rm * D^{k-i} \Rm + D^{k+3} T +
\sum_{i =
0}^{k} D^{i+1} T * D^{k-i} \Rm\\
&\ + \sum_{i = 0}^k \sum_{j = 0}^i D^j T * D^{i-j} T * D^{k-i} \Rm + \sum_{i =
0}^{k+2} D^i T * D^{k+2 - i} T.
\end{align*}
\begin{proof} We compute the first evolution equation, the case of $D^k T$ being
formally similar.  Using Proposition \ref{Nevolution} and Lemma \ref{generalS},
we compute
\begin{align*}
\dt D^k N =&\ \dt \left( \del + \gG \right) \dots \left(\del + \gG \right) N\\
=&\ \sum_{i = 0}^{k} D^i \left( \dt \gG \right)* D^{k-i} N + D^{k}
\left(\dt N \right)\\
=&\ \sum_{i = 0}^{k} D^i \left( \Rm + D T + T^{*2} \right) * D^{k-i} N + D^{k}
\left( \gD N + \Rm * T + D T * T + T^{*3} \right)\\
=&\ \gD D^k N + \sum_{i = 0}^{k-1} D^i \Rm * D^{k-i} T + \sum_{i = 1}^{k+1} D^i
T * D^{k+1-i} T + \sum_{i = 0}^{k} \sum_{j = 0}^i D^j T * D^{i-j}T * D^{k-i} T,
\end{align*}
as required.  Next we compute
\begin{align*}
\dt D^k \Rm =&\ \dt \left( \del + \gG \right) \dots \left(\del + \gG \right)
\Rm\\
=&\ \sum_{i = 0}^{k-1} D^i \left(\dt \gG \right) * D^{k-i - 1} \Rm + D^k
\left(\dt \Rm \right)\\
=&\ \sum_{i = 0}^{k-1} D^{i+1} \left(\Rm + D T + T^{*2} \right) * D^{k-i-1}
\Rm\\
&\ + D^k \left( \gD \Rm + \Rm^{*2} + \Rm * T^{*2} + \Rm * D T + T * D^2 T + D T
* D T + D^3 T \right)\\
=&\ \gD D^k \Rm + \sum_{i = 0}^k D^i \Rm * D^{k-i} \Rm + D^{k+3} T + \sum_{i =
0}^{k} D^{i+1} T * D^{k-i} \Rm\\
&\ + \sum_{i = 0}^k \sum_{j = 0}^i D^j T * D^{i-j} T * D^{k-i} \Rm + \sum_{i =
0}^{k+2} D^i T * D^{k+2 - i} T.
\end{align*}
\end{proof}
\end{thm}

%%%%%%%%%%%%%%%%%%%%%%%%%%%%%%%%%%%%%%%%%%%%%%%%%%%%%%%%%%%%%%
\subsection{Evolution equations for symplectic curvature flow}
%%%%%%%%%%%%%%%%%%%%%%%%%%%%%%%%%%%%%%%%%%%%%%%%%%%%%%%%%%%%%%

We begin by deriving an evolution equation for the Levi Civita derivative of
$J$.

\begin{prop} \label{sympDJev} Let $(M^{2n}, \omega(t), J(t))$ be a solution to
(\ref{AKflow}).  Then
\begin{gather} \label{sympDJevcalc}
\begin{split}
\dt \left(D J \right)_{ij}^k =&\ \gD D_i J_j^k {+ D_i \mathcal N_j^k} - g^{pl}
\left( D_i B_{jl} + D_j B_{il} - D_l B_{ij} \right) J_p^k\\
&\ + g^{kl} \left( D_i B_{pl} + D_p B_{il} - D_l B_{ip} \right) J_j^p\\
&\ + 2 g^{rs} \left( R_{s i j}^t D_r J_t^k - R_{s i t}^k D_r J_j^t \right) -
R_i^t D_t J_j^k + R_p^k D_i J_j^p - R_j^p D_i J_p^k,
\end{split}
\end{gather}
where
\begin{align*}
B = {\frac{1}{4} B^1 - \frac{1}{2} B^2}.
\end{align*}
\begin{proof} We recall the coordinate expression
\begin{align*}
D J_{ij}^k =&\ \del_i J_j^k - \gG_{i j}^p J_p^k + \gG_{i p}^k J_j^p
\end{align*}
If we choose normal coordinates for $g(0)$ at some point and differentiate this
expression with respect to $t$ this yields
\begin{align*}
\dt \left( DJ \right)_{ij}^k =&\ D_i \dot{J}_{j}^k - \dot{\gG}_{ij}^p J_p^k +
\dot{\gG}_{ip}^k J_j^p.
\end{align*}
Next we recall that if $\dt g = h$, one has
\begin{align*}
\dt \gG_{ij}^k =&\ \frac{1}{2} g^{kl} \left( D_i h_{jl} + D_j h_{il} - D_l
h_{ij} \right).
\end{align*}
To simplify the calculation, we set $B = {\frac{1}{4} B^1 - \frac{1}{2} B^2}$,
then using Proposition
\ref{Bfieldprop1} we derive
\begin{align*}
\dt \left( DJ \right)_{ij}^k =&\ D_i \left( g^{rs} D_r D_s J_j^k { + \mathcal
N_j^k} + \mathcal R_j^k \right)\\
&\ + g^{pl} \left( D_i \left(R_{jl} - B_{jl} \right) + D_j \left( R_{il} -
B_{il} \right) - D_l \left( R_{ij} - B_{ij} \right) \right) J_p^k\\
&\ - g^{kl} \left( D_i \left(R_{pl} - B_{pl} \right) + D_p \left( R_{il} -
B_{il} \right) - D_l \left( R_{ip} - B_{ip} \right) \right) J_j^p
\end{align*}
Now we commute derivatives and apply the differential Bianchi identity to
conclude
\begin{gather} \label{Bochner}
\begin{split}
D_i (g^{rs} D_r D_s J_j^k) =&\ g^{rs} \left(D_r D_i D_s J_j^k + R_{r i s}^t D_t
J_j^k + R_{r i j}^t D_s J_t^k - R_{ri t}^k D_s J_j^t \right)\\
=&\ g^{rs} \left( D_r \left(D_s D_i J_j^k + R_{s i j}^t J_t^k - R_{s i t}^k
J_j^t \right) \right)\\
&\ + g^{rs} \left(R_{r i s}^t D_t J_j^k + R_{r i j}^t D_s J_t^k - R_{ri t}^k
D_s J_j^t \right)\\
=&\ \gD D_i J_j^k + g^{rs} \left( D_r R_{s i j}^t J_t^k - D_r R_{s i t}^k J_j^t
\right)\\
&\ + g^{rs} \left( R_{s i j}^t D_r J_t^k - R_{s i t}^k D_r J_j^t + R_{r i s}^t
D_t J_j^k + R_{r i j}^t D_s J_t^k - R_{r i t}^k D_s J_j^t \right)\\
=&\ \gD D_i J_j^k + \left(D^t R_{i j} - D_j R_i^t \right) J_t^k - \left( D^k
R_{i t} - D_t R_i^k \right) J_j^t\\
&\ + 2 g^{rs} \left( R_{s i j}^t D_r J_t^k - R_{s i t}^k D_r J_j^t \right) -
R_i^t D_t J_j^k.
\end{split}
\end{gather}
Plugging this into the above line yields
\begin{align*}
\dt \left( DJ \right)_{ij}^k =&\ \gD D_i J_j^k  { + D_i \mathcal N_j^k} - g^{pl}
\left( D_i B_{jl} + D_j B_{il} - D_l B_{ij} \right) J_p^k\\
&\ + g^{kl} \left( D_i B_{pl} + D_p B_{il} - D_l B_{ip} \right) J_j^p\\
&\ + 2 g^{rs} \left( R_{s i j}^t D_r J_t^k - R_{s i t}^k D_r J_j^t \right) -
R_i^t D_t J_j^k + D_i \left( J_j^p R_p^k - R_j^p J_p^k \right)\\
&\ + \left(D^t R_{i j} - D_j R_i^t \right) J_t^k - \left( D^k R_{i t} - D_t
R_i^k \right) J_j^t\\
&\ + g^{pl} \left( D_i R_{jl} + D_j R_{il} - D_l R_{ij} \right) J_p^k\\
&\ - g^{kl} \left( D_i R_{pl} + D_p R_{il} - D_l R_{ip} \right) J_j^p\\
=&\ \gD D_i J_j^k {+ D_i \mathcal N_j^k} - g^{pl} \left( D_i B_{jl} + D_j B_{il}
-
D_l B_{ij} \right) J_p^k\\
&\ + g^{kl} \left( D_i B_{pl} + D_p B_{il} - D_l B_{ip} \right) J_j^p\\
&\ + 2 g^{rs} \left( R_{s i j}^t D_r J_t^k - R_{s i t}^k D_r J_j^t \right) -
R_i^t D_t J_j^k + R_p^k D_i J_j^p - R_j^p D_i J_p^k.
\end{align*}
Indeed, slightly miraculously, all of the $D R$ terms drop out of the equation.
This is the required result.
\end{proof}
\end{prop}

Next we use this proposition to compute the evolution of $\brs{DJ}^2$.

\begin{prop} \label{sympDJnormev}
Let $(M^{2n}, \omega(t), J(t))$ be a solution to
(\ref{AKflow}).  Then
\begin{gather} \label{sympDJcalc}
\begin{split}
\dt \brs{D J}^2 =&\ \gD \brs{D J}^2 - 2 \brs{D^2 J}^2 {- \frac{1}{2} \brs{B^1}^2
+ 3 \left< B^1, B^2 \right>}\\
&\ + 8 R_{s i j t} D_s J_t^k D_i J_j^k + {2 D_p D_i J_m^n D_k J_m^n J_j^p D_i
J_j^k + 2 D_p D_k J_m^n D_i J_m^n J_j^p D_i J_j^k}\\
&\ - {4 D_l D_m J_i^n D_m J_j^n J_l^k D_i J_j^k - 4 D_l D_m J_j^n D_m J_i^n
J_l^k
D_i J_j^k},
\end{split}
\end{gather}
where the result is interpreted in an orthonormal frame.
\begin{proof} First we note for $\dt g = h$ one has
\begin{align*}
\dt \brs{DJ}^2 =&\ \dt \left( g^{ip} g^{jq} g_{kr} D_i J_j^k D_p J_q^r \right)\\
=&\ 2 \left< \dt DJ, DJ \right> - \left< h, B^1 \right> + \left( g^{ip} g^{jq}
h_{kr} - g^{ip} h^{jq} g_{kr} \right) D_i J_j^k D_p J_q^r\\
=&\ 2 \left< \dt DJ, DJ \right> - \left< h, B^1 \right>.
\end{align*}
First of all, from (\ref{sympmetricev}), we have $h = -2 \Rc + \frac{1}{2} B^1 -
B^2$ and so
\begin{align} \label{DJnorm10}
- \left< h, B^1 \right> =&\ 2 \left< \Rc, B^1 \right> {- \frac{1}{2} \brs{B^1}^2
+ \left< B^1, B^2 \right>}.
\end{align}
We simplify the contributions of to $\dt DJ$ from (\ref{sympDJevcalc})
separately.  First note that
\begin{align} \label{DJnorm20}
2 \left< \gD D J, DJ \right> =&\ \gD \brs{DJ}^2 - 2 \brs{D^2 J}^2.
\end{align}
Next we simplify
\begin{align} \label{DJnorm30}
2 \left( -R_i^t D_t J_j^k + R_p^k D_i J_j^p - R_j^p D_i J_p^k \right) \left( D_i
J_j^k \right) =&\ -2 \left<\Rc, B^1 \right>.
\end{align}
Next we simplify, working in an orthonormal frame,
\begin{gather} \label{DJnorm40}
\begin{split}
4 g^{rs} \left(R_{s i j}^t D_r J_t^k - R_{s i t}^k D_r J_j^t \right) D_i J_j^k
=&\ 8 R_{s i j t} D_s J_t^k D_i J_j^k.
\end{split}
\end{gather}
Next consider
{
\begin{align*}
2 D_i \mathcal N_j^k D_i J_j^k =&\ 2 D_i \left(g^{kl} g_{mn} g^{pq} D_p
J_r^m J_j^r D_q J_l^n \right) D_i J_j^k\\
=&\ 2 \left[ D_i D_p J_r^m J_j^r D_p J_k^m + D_p J_r^m D_i J_j^r D_p J_k^m +
D_p J_r^m J_j^r D_i D_p J_k^m \right] D_i J_j^k.
\end{align*}
Combining the first and last terms and applying the identity $D (J^2) = 0$
yields
\begin{align} \label{magiccancel}
2 \left[ - J_j^k D_i D_p J_r^m D_i J_j^r D_p J_k^m + D_p J_r^m J_j^r D_i D_p
J_k^m D_i J_j^k\right] = 0.
\end{align}
Thus
\begin{gather} \label{DJnorm50}
2 D_i \mathcal N_j^k D_i J_j^k = 2 \left< B^1, B^2 \right>.
\end{gather}}
We simplify the remaining $DB$ terms of (\ref{sympDJevcalc}).  We begin by
considering the contribution of $B^1$.  First by a simplification like in line
(\ref{magiccancel}) we have
\begin{gather} \label{DJnorm60}
\begin{split}
- \frac{1}{2} D_i B^1_{jp} J_p^k D_i J_j^k =&\ - \frac{1}{2} D_i \left( D_j
J_m^n D_p J_m^n \right)
J_p^k D_i J_j^k\\
=&\ - \frac{1}{2} \left[D_i D_j J_m^n D_p J_m^n J_p^k D_i J_j^k + D_i D_p J_m^n
D_i J_p^k J_j^k
D_j J_m^n \right]\\
=&\ 0.
\end{split}
\end{gather}
Next we note
\begin{gather} \label{DJnorm70}
\begin{split}
- \frac{1}{2} D_j B^1_{ip} J_p^k D_i J_j^k =&\ - \frac{1}{2} D_j \left( D_i
J_m^n D_p J_m^n \right)
J_p^k D_i J_j^k\\
=&\ -  \frac{1}{2} D_j D_i J_m^n D_p J_m^n J_p^k D_i J_j^k - \frac{1}{2} D_j D_p
J_m^n D_i J_m^n J_p^k
D_i J_j^k.
\end{split}
\end{gather}
The next term is
\begin{gather} \label{DJnorm80}
\begin{split}
\frac{1}{2} D_l B^1_{ij} J_l^k D_i J_j^k =&\ \frac{1}{2} D_l \left( D_i J_m^n
D_j J_m^n \right) J_l^k
D_i J_j^k\\
=&\ \frac{1}{2} D_l D_i J_m^n D_j J_m^n J_l^k D_i J_j^k + \frac{1}{2} D_l D_j
J_m^n D_i J_m^n J_l^k D_i
J_j^k.
\end{split}
\end{gather}
Next
\begin{gather} \label{DJnorm90}
\begin{split}
\frac{1}{2} D_i B^1_{pk} J_j^p D_i J_j^k =&\ \frac{1}{2} D_i \left( D_p J_m^n
D_k J_m^n \right) J_j^p
D_i J_j^k\\
=&\ \frac{1}{2} D_i D_p J_m^n D_k J_m^n J_j^p D_i J_j^k + \frac{1}{2} D_i D_k
J_m^n D_p J_m^n J_j^p D_i
J_j^k\\
=&\ - \frac{1}{2} D_i D_p J_m^n D_k J_m^n D_i J_j^p J_j^k + \frac{1}{2} D_i D_k
J_m^n D_p J_m^n J_j^p
D_i J_j^k\\
=&\ 0.
\end{split}
\end{gather}
Next we compute
\begin{gather} \label{DJnorm100}
\begin{split}
\frac{1}{2} D_p B^1_{ik} J_j^p D_i J_j^k =&\ \frac{1}{2} D_p \left( D_i J_m^n
D_k J_m^n \right) J_j^p
D_i J_j^k\\
=&\ \frac{1}{2} D_p D_i J_m^n D_k J_m^n J_j^p D_i J_j^k + \frac{1}{2} D_p D_k
J_m^n D_i J_m^n J_j^p D_i
J_j^k.
\end{split}
\end{gather}
The final $B^1$ term is
\begin{gather} \label{DJnorm110}
\begin{split}
- \frac{1}{2} D_k B^1_{ip} J_j^p D_i J_j^k =&\ - \frac{1}{2} D_k \left( D_i
J_m^n D_p J_m^n \right)
J_j^p D_i J_j^k\\
=&\ - \frac{1}{2} D_k D_i J_m^n D_p J_m^n J_j^p D_i J_j^k - \frac{1}{2} D_k D_p
J_m^n D_i J_m^n J_j^p
D_i J_j^k.
\end{split}
\end{gather}
Now we compute the $B^2$ terms.  The first such is
\begin{gather} \label{DJnorm120}
\begin{split}
D_i B^2_{jp} J_p^k D_i J_j^k =&\ D_i \left( D_m J_j^n D_m J_p^n \right)
J_p^k D_i J_j^k\\
=&\ D_i D_m J_j^n D_m J_p^n J_p^k D_i J_j^k + D_i D_m J_p^n D_m J_j^n J_p^k
D_i J_j^k\\
=&\ D_i D_m J_j^n D_m J_p^n J_p^k D_i J_j^k - D_i D_m J_p^n D_m J_j^n D_i
J_p^k J_j^k\\
=&\ 0.
\end{split}
\end{gather}
Next we note
\begin{gather} \label{DJnorm130}
\begin{split}
 D_j B^2_{ip} J_p^k D_i J_j^k =&\ D_j \left( D_m J_i^n D_m J_p^n \right)
J_p^k D_i J_j^k\\
=&\ D_j D_m J_i^n D_m J_p^n J_p^k D_i J_j^k + D_j D_m J_p^n D_m J_i^n J_p^k
D_i J_j^k.
\end{split}
\end{gather}
The next term is
\begin{gather} \label{DJnorm140}
\begin{split}
- D_l B^2_{ij} J_l^k D_i J_j^k =&\ - D_l \left( D_m J_i^n D_m J_j^n \right)
J_l^k
D_i J_j^k\\
=&\ - D_l D_m J_i^n D_m J_j^n J_l^k D_i J_j^k - D_l D_m J_j^n D_m J_i^n J_l^k
D_i
J_j^k.
\end{split}
\end{gather}
Next
\begin{gather} \label{DJnorm150}
\begin{split}
- D_i B^2_{pk} J_j^p D_i J_j^k =&\ - D_i \left( D_m J_p^n D_m J_k^n \right)
J_j^p
D_i J_j^k\\
=&\ - D_i D_m J_p^n D_m J_k^n J_j^p D_i J_j^k - D_i D_m J_k^n D_m J_p^n J_j^p
D_i
J_j^k\\
=&\ - D_i D_m J_p^n D_m J_k^n J_j^p D_i J_j^k + D_i D_m J_k^n D_m J_p^n J_j^k
D_i
J_j^p\\
=&\ 0
\end{split}
\end{gather}
Next we compute
\begin{gather} \label{DJnorm160}
\begin{split}
- D_p B^2_{ik} J_j^p D_i J_j^k =&\ - D_p \left( D_m J_i^n D_m J_k^n \right)
J_j^p
D_i J_j^k\\
=&\ - D_p D_m J_i^n D_m J_k^n J_j^p D_i J_j^k - D_p D_m J_k^n D_m J_i^n J_j^p
D_i
J_j^k.
\end{split}
\end{gather}
The final $B^2$ term is
\begin{gather} \label{DJnorm170}
\begin{split}
D_k B^2_{ip} J_j^p D_i J_j^k =&\ D_k \left( D_m J_i^n D_m J_p^n \right)
J_j^p D_i J_j^k\\
=&\ D_k D_m J_i^n D_m J_p^n J_j^p D_i J_j^k + D_k D_m J_p^n D_m J_i^n J_j^p
D_i J_j^k.
\end{split}
\end{gather}
By further applying the identity $D J^2 = 0$ one may observe that
lines (\ref{DJnorm70}), (\ref{DJnorm80}), (\ref{DJnorm100}), (\ref{DJnorm110})
are all equal.  Likewise the lines (\ref{DJnorm130}), (\ref{DJnorm140}),
(\ref{DJnorm160}), (\ref{DJnorm170}) are equal.  Combining the calculations
yields the result.
\end{proof}
\end{prop}

Next we derive an evolution equation for the Riemannian curvature tensor.
\begin{prop} \label{sympcurvev} Let $(M^{2n}, \omega(t), J(t))$ be a solution to
(\ref{AKflow}).  Then
\begin{gather} \label{sympcurvevcalc}
\begin{split}
\dt R_{ijkl} =&\ \gD R_{ijkl} + \Rm_{ijkl}^2 + \Rm_{ijkl}^{\sharp} -
\left(R_{ip}R_{pjkl} + R_{jp} R_{ipkl} + R_{kp} R_{ijpl} + R_{lp} R_{ijkp}
\right)\\
&\ + D_i D_k B_{jl} - D_i D_l B_{jk} + D_j
D_k B_{il} + D_j D_l B_{ik} + B_{l}^p R_{ijkp} - B_k^p R_{ijlp},
\end{split}
\end{gather}
where
\begin{align*}
B = {\frac{1}{4} B^1 - \frac{1}{2}  B^2}.
\end{align*}
\begin{proof} The starting point is the formula (\ref{curvvar}).  A well-known
calculation (\cite{Chow} Lemma 2.51) shows that for $h = -2 \Rc$,
\begin{align*}
\dt \Rm_{ijkl} =&\ \gD R_{ijkl} + \Rm_{ijkl}^2 + \Rm_{ijkl}^{\sharp} -
\left(R_{ip}R_{pjkl} + R_{jp} R_{ipkl} + R_{kp} R_{ijpl} + R_{lp} R_{ijkp}
\right)
\end{align*}
where $\Rm^2$ is the square of the curvature operator and $\Rm^{\sharp}$ is the
Lie algebra square.  Including the extra term of $h = 2 B$, the result follows.
\end{proof}
\end{prop}

In the next theorem we derive evolution equations for the covariant derivatives
of the Nijenhuis tensor and curvature for solutions to (\ref{AKflow}).
\begin{thm} \label{sympderev} Let $(M^{2n}, \omega(t), J(t))$ be a solution to
(\ref{AKflow}).  Then
\begin{align*}
\dt D^k J =&\ \gD D^k J + \sum_{i = 1}^{k} D^i J * D^{k-i} \Rm + \sum_{i =
1}^{k+1} D^i J * D^{k+2 - i} J\\
&\ + \sum_{i = 0}^{k-1} \sum_{j = 0}^i D^{j+1} J *D^{i+1-j} J* D^{k-i} J,\\
\dt D^k \Rm =&\ \gD D^k \Rm + \sum_{i = 0}^k D^i \Rm * D^{k-i} \Rm\\
&\ + \sum_{i=0}^k \sum_{j=0}^i D^{k-i} \Rm* D^{j+1} J *D^{i-j+1} J + \sum_{i =
0}^k D^{i+1} J * D^{k-i+3} J.
\end{align*}
\begin{proof} We express $D^{k} J = D^{k-1} D J$ and compute using the result of
Proposition \ref{sympDJev}
\begin{align*}
\dt D^k J =&\ \dt \left( \del + \gG \right) \dots \left(\del + \gG \right) D J\\
=&\ \sum_{i = 0}^{k-2} D^i \left( \dt \gG \right) * D^{k-i-2} DJ + D^{k-1}
\left(\dt DJ \right)\\
=&\ \sum_{i = 0}^{k-2} D^{i+1} \left( \Rm + DJ^{*2} \right) * D^{k-i-1} J +
D^{k-1} \left( \gD DJ + D^2 J * DJ + \Rm * DJ \right)\\
=&\ \gD D^k J + \sum_{i = 1}^{k} D^i J * D^{k-i} \Rm + \sum_{i = 0}^{k-1}
\sum_{j = 0}^i D^{j+1} J D^{i+1-j} J D^{k-i} J\\
&\ + \sum_{i = 1}^{k+1} D^i J * D^{k+2 - i} J.
\end{align*}
Likewise we compute using the result of Proposition \ref{sympcurvev}),
\begin{align*}
\dt D^k \Rm =&\ \dt \left( \del + \gG \right) \dots \left(\del + \gG \right)
\Rm\\
=&\ \sum_{i = 0}^{k-1} D^i \left( \dt \gG \right)* D^{k-i-1} \Rm + D^k \left(
\dt
\Rm \right)\\
=&\ \sum_{i = 0}^{k-1} D^{i+1} \left( \Rm + DJ^{*2} \right)* D^{k-i-1} \Rm\\
&\ + D^k \left( \gD \Rm + \Rm^{*2} + \Rm * DJ^{*2} + D^2 J * D^2 J + D^3 J * D J
\right)\\
=&\ \gD D^k \Rm + \sum_{i = 0}^k D^i \Rm * D^{k-i} \Rm\\
&\ + \sum_{i=0}^k \sum_{j=0}^i D^{k-i} \Rm *D^{j+1} J *D^{i-j+1} J + \sum_{i =
0}^k D^{i+1} J * D^{k-i+3} J,
\end{align*}
as required.
\end{proof}
\end{thm}

%%%%%%%%%%%%%%%%%%%%%%%%%%%%%%%%%%%%%%%%%%%%%%%%%%%%%%%%%%%%%%%%%%%
\section{Smoothing Estimates and a Long time existence obstruction}
\label{smoothing}
%%%%%%%%%%%%%%%%%%%%%%%%%%%%%%%%%%%%%%%%%%%%%%%%%%%%%%%%%%%%%%%%%%%

We begin with the proof of Theorem \ref{generalsmoothingthm}.

\begin{proof} Arguments of this kind are standard.  We give the proof for the
case $m = 1$ to indicate why we are forced to deal with $L^2$ norms for this
flow.  Let $A$ denote a constant to be determined later, and consider
\begin{align*}
F(t) := t \left( \brs{\brs{D \Rm}}_{L^2}^2 + \brs{\brs{D^2 T}}_{L^2}^2 \right) +
\ga \brs{\brs{\Rm}}_{L^2}^2 + \gb \brs{\brs{D T}}_{L^2}^2 + \gamma \brs{\brs{D
N}}_{L^2}^2.
\end{align*}
Using the result of Theorem \ref{ACFderevolution}, integrating by parts and
repeatedly applying the Cauchy-Schwarz inequality yields
\begin{align*}
\frac{d}{dt} F \leq&\ \left(\brs{\brs{D \Rm}}_{L^2}^2 + \brs{\brs{D^2
T}}_{L^2}^2 \right) \left(1 + CtK \right)\\
&\ - 2 \ga \brs{\brs{ D \Rm}}_{L^2}^2 -2 \gb \brs{\brs{D^2 T}}_{L^2}^2 - 2
\gamma \brs{\brs{D^2 N}}_{L^2}^2\\
&\ + \ga \int_M D^3 T * \Rm + \gb \int_M D^3 N * DT + C \left(1 + tK \right) F.
\end{align*}
The first two terms in the last line are reflecting the upper triangularity of
the symbol
of our system of equations, and are the reason for resorting to $L^2$ norms.
Integrating by parts and applying the Cauchy Schwarz inequality, we observe that
with $\ga$ chosen large with respect to a uniform constant, $\gb$ chosen large
with respect to $\ga$, and $\gamma$ chosen large with respect to $\gb$ we yield
the estimate
\begin{align*}
\frac{d}{dt} F \leq&\ C(\left(1 + t K \right)) F
\end{align*}
which one may integrate over the required time interval to yield the required
result.
\end{proof}

\noindent Next we claim pointwise smoothing estimates for solutions to
(\ref{AKflow}).

\begin{thm} \label{AKsmoothingthm} Given $m > 0$, there exists $C = C(m,
n)$ such that if $(M^{2n}, \omega(t), J(t))$ is a solution to
(\ref{AKflow}) on $\left[0, \frac{\ga}{K} \right]$ satisfying
\begin{align*}
\sup_{M \times \left[0, \frac{\ga}{K} \right]} \{ \brs{\Rm}, \brs{DJ}^2,
\brs{D^2
T}
\} \leq K,
\end{align*}
then
\begin{align*}
\sup_{M \times \left(0, \frac{\ga}{K} \right]} \{ \brs{ D^m \Rm}_{C^0},
\brs{D^{m+2} J}_{C^0} \} \leq \frac{C K}{t^{\frac{m}{2}}}.
\end{align*}
\begin{proof} Though the quantities involved are different, the proof is
formally identical to \cite{ST1} Theorem 7.3, with $DJ$ playing the role of $T$.
\end{proof}
\end{thm}

\noindent We can now give the proof of Theorem \ref{lteobs}
\begin{proof} This argument is also standard.  First one notes that if the
curvature, torsion, and first derivative of torsion are bounded on a finite time
interval, then the metrics are uniformly equivalent along the flow, and hence
the Sobolev constant of the manifold is also bounded.  The $L^2$ derivative
bounds of Theorem \ref{ACFderevolution} then yield pointwise bounds on the
derivatives of curvature and torsion.  Then an application of
\cite{HamiltonComp} Theorem 2.3 yields smooth convergence of the almost
Hermitian
structures as $t \to \tau$, contradicting maximality of $\tau$.
\end{proof}

Next we want to prove Theorem \ref{symplte}, which will follow directly
from Theorem \ref{lteobs} and the following proposition.
{
\begin{prop} \label{torsionbound} Let $(M^{2n}, \omega(t), J(t))$ be a solution
to
(\ref{AKflow}) on
$[0, T]$ satisfying
\begin{align*}
\sup_{M \times [0, T]} \brs{\Rm} = K.
\end{align*}
There exists a constant $C(K, n, \omega_0, J_0,T)$ such that
\begin{align*}
\sup_{M \times [0, T]} \brs{DJ}^2 + \brs{D^2 J} \leq C.
\end{align*}
\begin{proof} From the Weitzenb\"ock formula (Lemma \ref{weitzenbock}) we note
that for an almost K\"ahler manifold with $\brs{\Rm} \leq K$, one has $\brs{D^*
D \omega} \leq C K$ for a constant $C$ depending only on $n$.  In the proof of
Lemma \ref{P2formula} we noted that $\left( D^* D \omega \right)^{1,1} = N^2$. 
By a direct calculation one has that $\left< \omega, N^2 \right> = \tr_g B^2 =
\brs{DJ}^2$.  Thus, one has
\begin{align*}
\brs{DJ}^2 =\brs{\left<\omega, N^2 \right>} =&\ \brs{\left< \omega, D^*D \omega
\right>} \leq C \brs{D^* D \omega} \leq C K.
\end{align*}
Now fix some constants $\ga, \gb > 0$ and let
\begin{align*}
F(x, t) =&\ \brs{D^2 J}^2 + \ga \brs{\Rm}^2 + \gb \brs{DJ}^2.
\end{align*}
Using Proposition \ref{sympcurvev} and Theorem \ref{sympderev} we conclude that
there is a constant $C$ depending only on $n$ such that
\begin{align*}
\dt F \leq&\ \gD F - 2 \left( \brs{D^3 J}^2 + \ga \brs{D \Rm}^2 + \gb \brs{D^2
J}^2 \right)\\
&\ + C \left( \brs{\Rm} \brs{D^2 J}^2 + \brs{D J} \brs{D \Rm} \brs{D^2 J} +
\brs{D J} \brs{D^2 J} \brs{D^3 J} + \brs{D J}^2 \brs{D^2 J}^2 \right)\\
&\ + C \ga \left( \brs{\Rm}^3 + \brs{D J} \brs{D^3 J} \brs{\Rm} + \brs{D^2 J}^2
\brs{\Rm} + \brs{DJ}^2 \brs{\Rm}^2 \right)\\
&\ + C \gb \left( \brs{DJ}^4 + \brs{\Rm} \brs{DJ}^2 + \brs{D^2 J} \brs{DJ}^2
\right).
\end{align*}
Using the estimates on $\brs{\Rm}$ and $\brs{DJ}^2$, and replacing $K$ by $\max
\{K, 1\}$ if necessary we conclude by the Cauchy-Schwarz inequality
\begin{align*}
\dt F \leq&\ \gD F- 2 \left( \brs{D^3 J}^2 + \ga \brs{D \Rm}^2 + \gb \brs{D^2
J}^2 \right)\\
&\ + \frac{1}{2} \brs{D^3 J}^2 + C K \brs{D^2 J}^2 + C K^2 \left( \brs{D \Rm}^2
+ \brs{D^2 J}^2 \right)\\
&\ + \frac{1}{2} \brs{D^3 J}^2 + C \ga K \brs{D^2 J}^2 + C \ga^2 K^3\\
&\ + \frac{1}{2} \brs{D^2 J}^2 + C \left(\gb + \gb^2 \right) K^2.
\end{align*}
It is clear then that if we choose $\ga$ large with respect to dimension
dependent constants, and then choose $\gb$ large with respect to $\ga$, we may
conclude
\begin{align*}
\dt F \leq&\ \gD F + C K^3.
\end{align*}
The proposition follows from the maximum principle.
\end{proof}
\end{prop}}

%%%%%%%%%%%%%%%%%%%%%%%%%%%%%%%%%%%%%%%%%%%
\section{The structure of critical metrics} \label{critstruct}
%%%%%%%%%%%%%%%%%%%%%%%%%%%%%%%%%%%%%%%%%%%

In this section we record some results on the structure of the limiting objects
of equations (\ref{AKflow}).

\begin{defn} Let $(M^{2n}, \omega, J)$ be an almost K\"ahler manifold.  We say
that this manifold is \emph{static} if there exists $\gl \in \mathbb R$ such
that
\begin{align} \label{static1}
P =&\ \gl \omega\\
D^* D J - \mathcal N - \mathcal R =&\ 0. \label{static2}
\end{align}
\end{defn}

Let us say a word on the definition of this condition.  We want to understand
the limiting behavior of equation (\ref{AKflow}), hence the first condition
arises for solutions which simply rescale the metric.  Observe though that even
for solutions which are scaling the metric, one expects $J$ to remain fixed as
one cannot scale almost complex structures.  Thus the static condition defined
above is a natural expression of the expected smooth limit points of
(\ref{AKflow}).

\begin{lemma} \label{Jsym} Let $(M^{2n}, \omega, J)$ be a static structure. Then
\begin{align} \label{RicJinv}
\Ric^{-J} =&\ 0,
\end{align}
i.e. the Ricci tensor is $J$-invariant.
\begin{proof} Equation (\ref{static1}) implies that $P^{2,0 + 0,2} = 0$.
Equation (\ref{static2}) may be expressed as
\begin{align*}
g^{-1} \left[ P^{2,0 + 0,2} + \Ric^{-J} \right] = 0,
\end{align*}
and so the lemma follows.
\end{proof}
\end{lemma}

Let us show some further structure in dimension $4$.  Let $(M^4, g)$ be an
oriented Riemannian manifold.  Since one may decompose $\Lambda^2 = \Lambda^+
\oplus \Lambda^-$, the action of the curvature tensor on $\Lambda^2$ decomposes
accordingly, and is typically written
\begin{align} \label{basicdecomp}
R = \left(
\begin{array}{c|c}
W^+ + \frac{s}{12} I & \ohat{\Rc}\\
\hline
\ohat{\Rc} & W^- + \frac{s}{12} I
\end{array} \right)
\end{align}
where $\ohat{\Rc}$ is a certain action of the traceless Ricci tensor and $W^+$
and $W^-$ are the self-dual and anti-self-dual Weyl curvatures.  If one further
has $(M^4, \omega, J)$ an almost Hermitian manifold, then one can refine the
decomposition of $\Lambda^2$ as
\begin{align} \label{refinedforms}
\Lambda^2 = \left( (\omega) \oplus \Lambda^{2,0} \right) \oplus \Lambda_0^{1,1}
\end{align}
where $\Lambda_0^{1,1}$ are real $(1,1)$ forms orthogonal to $\omega$.  Using
this further decomposition one yields, adopting notation of \cite{Armstrong},
\begin{align} \label{refineddecomp}
R = \left(
\begin{array}{c|c||c}
a & W^+_F & R_F\\
\hline
W_F^{+*} & W_{00}^+ + \frac{1}{2} b I & R_{00}\\
\hline
\hline
R_F^* & R_{00}^* & W_{00}^- + \frac{1}{3} c I
\end{array} \right)
\end{align}
where the tensors in this equation are defined by comparing with
(\ref{basicdecomp}) and using the refined decomposition of forms of
(\ref{refinedforms}).  The double
bars indicate the original decomposition into self-dual and anti-self-dual
forms.  Now we recall a curvature calculation in \cite{Armstrong} which
decomposes the
curvature tensor of the \emph{canonical} connection of an almost K\"ahler
manifold
according to (\ref{refinedforms}).

\begin{prop} \label{armstrong} (\cite{Armstrong} Proposition 2)
\begin{align} \label{Omegadecomp}
\Omega = \left(
\begin{array}{c|c||c}
\frac{s^{\N}}{12} & W^+_F & R_F - 2 C\\
\hline
0 & 0 & 0\\
\hline
\hline
R_F^* & R_{00} & W^- + \frac{1}{3} c I
\end{array} \right)
\end{align}
\end{prop}
\noindent One may consult \cite{Armstrong} for the precise definition of $C$,
which is not relevant to us here.  All the other tensors are the same as what
appears in (\ref{refineddecomp}).  It is important to observe that this matrix
acts from the right on two-forms.  For instance, the image acting from the right
lies entirely in $(1,1)$ forms, as required.

\begin{prop} \label{staticEin} Let $(M^4, \omega, J)$ be a static structure.
Then
\begin{align} \label{weylflat}
W_F \equiv 0.
\end{align}
\begin{proof} This immediate from (\ref{Omegadecomp}) and the fact that $P = \gl
\omega$.
\end{proof}
\end{prop}

Returning to (\ref{refineddecomp}) it follows that the $\omega$ is an
eigenvector for the action of $W_+$.  This condition is related to delicate
topological estimates of LeBrun \cite{LeBrun} related to the Seiberg-Witten
equations.  Furthermore, a Theorem of Apostolov, Armstrong, and Dr\`aghici
states that compact almost K\"ahler four-manifolds satisfying (\ref{RicJinv})
and (\ref{weylflat}) are automatically K\"ahler.  Thus using Lemma
\ref{Jsym} and Proposition \ref{staticEin}, and \cite{Apostolov2} Theorem 2 it
follows that
compact static four-manifolds are K\"ahler Einstein.
\begin{cor} Let $(M^4, \omega, J)$ be a compact static structure.  Then
$(\omega, J)$ is K\"ahler-Einstein.
\end{cor}

%%%%%%%%%%%%%%%%%%%%%%%%%%%%%%%%%%%
\section{Remarks and open problems} \label{conc}
%%%%%%%%%%%%%%%%%%%%%%%%%%%%%%%%%%%

Recall from \cite{Tian07}, \cite{TZ} we know that the solution to K\"ahler Ricci
flow exists smoothly as long as the associated cohomology class is in the
K\"ahler cone.  Therefore it is natural, for purposes of understanding the long
time existence and singularity formation of solutions to (\ref{AKODE}), to
understand the corresponding cone $\mathcal C$ of symplectic forms in $H^2(M,
\mathbb R)$. Note that $\mathcal C$ consists of all cohomology classes
in $H^2(M, \mathbb R)$ which can be represented by a symplectic form.
Any symplectic form $\omega$ admits compatible almost complex structures, and
moreover the
space of these almost complex structures is contractible.  Thus one may define
the \emph{canonical class}
\begin{align*}
K = c_1(M, \omega) := c_1(M, J)
\end{align*}
where $J$ is any almost complex structure compatible with $\omega$ and the
orientation. It is clear that the homotopy classes of
symplectic structures define the same canonical class.  Therefore, associated to
a solution of (\ref{AKflow}), one has the well-defined \emph{associated ODE in
cohomology}
\begin{align} \label{AKODE}
\frac{d}{dt} [\omega] = - K.
\end{align}
It is clear by the definition that given a solution to (\ref{AKflow}), the
associated one parameter family of cohomology classes satisfies (\ref{AKODE}).
Thus, we have
\begin{lemma} Given $(M^{2n}, \omega(t), J(t))$ a solution to (\ref{AKflow}),
let
\begin{align*}
T^* := \sup \{ t > 0 | [\omega(t)] = [\omega(0)] - t K \in \mathcal C \}.
\end{align*}
Furthermore, let $T$ denote the maximal existence time of $(\omega(t), J(t))$.
Then
\begin{align*}
T \leq T^*
\end{align*}
\end{lemma}
\noindent
It is natural to conjecture: {\it The maximal existence time for (\ref{AKflow})
with initial condition $\omega(0)$ is given by $T^*$}.
This is the
analogue of the theorem of Tian-Zhang (\cite{Tian07}, \cite{TZ}) mentioned above
for K\"ahler Ricci flow.

If the above $T^* < \infty$, then (\ref{AKflow}) develops finite-time
singularity. The second basic problem is to study
{\it the nature of such a singularity}. Is it possible that such a singularity
is caused by $J$-holomorphic spheres as we see
in the case of K\"ahler manifolds? The case of 4-dimensional symplectic
manifolds is of particular interest and may be easier to study. We expect that
{\it either $(\omega(t), J(t))$ collapses to a lower dimensional space or
converges to a smooth pair $(\omega_T, J_T)$ outside a subvariety as $t $ tends
to $T\leq T^*$}. If so, we may do surgery and extend
(\ref{AKflow}) across $T$.  By scaling, one may get ancient solutions for
(\ref{AKflow}). A
basic problem is to classify all the ancient solutions.  In dimension 4, it may
be
possible to classify.

Another natural problem is to find functionals which are monotonic along
(\ref{AKflow}).   In particular one can ask, {\it is (\ref{AKflow}) a gradient
flow
like the Ricci flow and the pluriclosed flow of }\cite{ST2}?  We showed in
\cite{ST3} that the parabolic flow of pluriclosed metrics of \cite{ST2} is in
fact a gradient flow.  This was done by exhibiting that after change by a
certain diffeomorphism solutions to this flow are equivalent to solutions to the
$B$-field renormalization group flow of string theory.  In light of Proposition
\ref{Bfieldprop1}, solutions to (\ref{AKflow}) have the metric evolving by the
Ricci flow plus certain lower order terms, therefore one expects to be able to
add a certain Lagrangian to the Perelman functionals to obtain a gradient flow
property for (\ref{AKflow}), as in the $B$-field renormalization group flow.
This will be the subject of future work.

Furthermore, we believe that this new symplectic curvature flow will be useful
in
studying the topology of symplectic manifolds, particularly in dimension 4.
It follows from the results in section 6 that static
solutions in dimension 4 are of anti-self-dual type, more precisely, the
self-dual
part of curvature for the canonical connection is determined by its scalar
curvature. This gives a hope to use (\ref{AKflow}) to prove a symplectic version
of the Miyaoka-Yau inequality for complex surfaces. Such an inequality for
symplectic 4-manifolds has been long speculated.  For still further
applications, we are led to studying limits of (\ref{AKflow}) as time $t$ tends
to $\infty$ and after appropriate scalings. The limits should include the above
static metrics, soliton solutions as well as
collapsed metrics which generalize the metrics studied by Song-Tian
\cite{SongTian} for elliptic
surfaces.

Finally, one could go still further and attempt to use symplectic curvature flow
to understand general four-manifolds with $b_2^+ \geq 1$.  A smooth $4$-manifold
$M$ with $b_2^+ \geq 1$, admits a ``near-symplectic form,'' roughly speaking a
symplectic form which degenerates along certain disjoint circles in the
manifold.  Forms
of this type have been studied for instance in \cite{Taubes3}.  By imposing the
appropriate boundary condition, one may be able to construct solutions of
(\ref{AKflow}) for such near-symplectic structures.  In this manner one
envisions using symplectic curvature flow as a tool for understanding more
general $4$-manifolds.

\bibliographystyle{hamsplain}
%%%%%%%%%%%%%%%%%%%%%%%%%%%

\end{document}